\documentclass[10pt]{amsart}
\usepackage{amscd}
\usepackage{amssymb}
\date{14 July 2004}
\title[Differential Algebras]{Dualizing Complexes and 
Perverse Modules over Differential Algebras}

\author{Amnon Yekutieli and James J. Zhang }
\address{A. Yekutieli: Department of  Mathematics 
Ben Gurion University, 
Be'er Sheva 84105, 
Israel}
\email{amyekut@math.bgu.ac.il}
\address{J.J. Zhang: Department of Mathematics, Box 354350,
University of Washington, Seattle, Washington 98195, USA}
\email{zhang@math.washington.edu}
\thanks{{\em Mathematics Subject Classification} 2000.
Primary: 16D90; Secondary: 18G10, 16S32, 16W70, 16U20.}
\keywords{Noncommutative rings, filtered rings,
dualizing complexes.}
\thanks{This research was supported by the US-Israel Binational
Science Foundation. The second author was partially supported by 
the US National Science Foundation.}

\newtheorem{thm}[equation]{Theorem}
\newtheorem{cor}[equation]{Corollary}
\newtheorem{prop}[equation]{Proposition}
\newtheorem{lem}[equation]{Lemma}
\theoremstyle{definition}
\newtheorem{dfn}[equation]{Definition}
\newtheorem{rem}[equation]{Remark}
\newtheorem{exa}[equation]{Example}

\newtheorem{que}[equation]{Question}

\numberwithin{equation}{section}

\newcommand{\iso}{\xrightarrow{\simeq}}
\newcommand{\inj}{\hookrightarrow}
\newcommand{\surj}{\twoheadrightarrow}
\newcommand{\xar}{\xrightarrow}
\newcommand{\opn}{\operatorname}
\newcommand{\opnt}[1]{\mathrm{#1}} 
\newcommand{\cat}[1]{\operatorname{\mathsf{#1}}}

\newcommand{\rmitem}[1]{\item[\text{\textup{(#1)}}]}
\newcommand{\mfrak}[1]{\mathfrak{#1}}
\newcommand{\mcal}[1]{\mathcal{#1}}
\newcommand{\msf}[1]{\mathsf{#1}}

\newcommand{\mrm}[1]{\mathrm{#1}}
\newcommand{\mbb}[1]{\mathbb{#1}}

\newcommand{\tup}[1]{\textup{#1}}
\newcommand{\bsym}[1]{\boldsymbol{#1}}
\newcommand{\boplus}{\bigoplus\nolimits}
\newcommand{\wtil}[1]{\widetilde{#1}}

\newcommand{\bra}[1]{\langle #1 \rangle}
\renewcommand{\k}{\Bbbk}
\newcommand{\bwedge}{\bigwedge\nolimits}

\begin{document}

\begin{abstract}
A differential algebra of finite type over a field $\k$ is 
a filtered algebra $A$, such that the associated graded 
algebra is finite over its center, and the center is  
a finitely generated $\k$-algebra. The prototypical 
example is the algebra of differential operators on a 
smooth affine variety, when $\opn{char} \k = 0$. We 
study homological and geometric properties of 
differential algebras of finite type. The main 
results concern the rigid dualizing complex over such an
algebra $A$: its existence, structure and variance properties. 
We also define and study perverse $A$-modules,
and show how they are related to the Auslander 
property of the rigid dualizing complex of $A$. 
\end{abstract}

\maketitle

\setcounter{section}{-1}
\section{Introduction}

The ``classical'' Grothendieck Duality theory, dealing with 
dualizing complexes over schemes, was developed in 
the book {\em Residues and Duality} by Hartshorne \cite{RD}. 
A duality theory for noncommutative noetherian algebras over 
a base field $\k$ was introduced in \cite{Ye1}. Roughly 
speaking, a dualizing complex over a $\k$-algebra $A$ is a complex
$R \in \msf{D}^{\mrm{b}}(\cat{Mod} A^{\mrm{e}})$, 
such the functor 
\[ \mrm{RHom}_{A}(-, R) : \msf{D}^{\mrm{b}}_{\mrm{f}}(\cat{Mod} A)
\to \msf{D}^{\mrm{b}}_{\mrm{f}}(\cat{Mod} A^{\mrm{op}}) \]
is a duality, with inverse 
$\mrm{RHom}_{A^{\mrm{op}}}(-, R)$.
Here $A^{\mrm{op}}$ is the opposite algebra,
$A^{\mrm{e}} := A \otimes_{\k} A^{\mrm{op}}$,
and $\msf{D}^{\mrm{b}}_{\mrm{f}}(\cat{Mod} A)$
is the derived category of bounded complexes of $A$-modules 
with finite cohomologies. See Definition \ref{dfn2.1} for details.
In the decade since its introduction this noncommutative duality
theory has progressed in several directions; cf.\ 
the papers \cite{VdB}, \cite{Jo}, \cite{MY}, \cite{WZ}
and \cite{Ch2}. 

One of the biggest problems in noncommutative duality theory 
is {\em existence of dualizing complexes}. 
The most effective existence criterion to date is 
due to Van den Bergh \cite{VdB}. It says the following:
suppose the $\k$-algebra $A$ has some exhaustive nonnegative 
filtration $G = \{ G_i A \}$ such that the graded algebra 
$\opnt{gr}^G A$ is connected, noetherian and commutative. 
Then $A$ has a rigid dualizing complex $R_A$ 
(see Definition \ref{dfn2.3}). Moreover in \cite{YZ1} we proved that
the dualizing complex $R_A$ has the Auslander property (see 
Definition \ref{dfn2.2}), and it is unique up to a unique 
isomorphism.

Let us remind the reader that a graded $\k$-algebra $\bar{A}$ is 
called  {\em connected} if 
$\bar{A} = \bigoplus_{i \in \mbb{N}} \bar{A}_i$, 
$\bar{A}_0 = \k$ and each $\bar{A}_i$ is a finite $\k$-module. 
An exhaustive nonnegative  
filtration $G$ on $A$, such that the graded algebra 
$\bar{A} := \opnt{gr}^G A$ 
is a noetherian connected graded $\k$-algebra, can 
be considered as a ``noncommutative compactification of
$\opn{Spec} A$''. Indeed, when $A$ is commutative, let 
$\wtil{A} := \opnt{Rees}^G A \subset A[t]$ be the Rees algebra. 
Then $\opn{Proj} \wtil{A}$ is a projective scheme over $\k$, the 
divisor $\{ t = 0 \}$ is ample, and its complement is isomorphic 
to $\opn{Spec} A$.

We have observed that often in the literature one encounters 
algebras $A$ that are equipped with an exhaustive 
nonnegative filtration $F = \{ F_i A \}$, 
such that $\opnt{gr}^F A$ is noetherian and
finite over its center, yet is not connected. 
The main goal of this paper is to prove that Van den Bergh's 
existence criterion applies to such algebras too, 
and furthermore the rigid dualizing complex $R_A$ has 
especially good homological and geometric properties. These 
properties shall be  used in our sequel paper \cite{YZ3} to 
construct rigid dualizing complexes over noncommutative ringed 
schemes.

Let us introduce some conventions now. 
Throughout the paper $\k$ is a field. By default all $\k$-algebras 
are associative unital algebras, all modules are left modules,
and all bimodules are central over $\k$. Given a $\k$-algebra $A$
we denote by $\cat{Mod} A$ the category of 
$A$-modules.  The unadorned tensor product $\otimes$ will mean 
$\otimes_{\k}$. 

Let $C$ be a finitely generated 
commutative $\k$-algebra and let $A$ be
a $C$-ring (i.e.\ there is a ring homomorphism $C \to A$). 
We call $A$ a {\em differential $C$-ring of finite type} if there 
exists a nonnegative exhaustive filtration 
$F = \{ F_i A \}_{i \in \mbb{Z}}$
of $A$ such that $\opnt{gr}^{F} A$ is a finite 
module over its center $\mrm{Z}(\opnt{gr}^{F} A)$, and 
$\mrm{Z}(\opnt{gr}^{F} A)$ is a finitely generated $C$-algebra. 
We also call $A$ a {\em differential $\k$-algebra of finite type}. 
See Definitions \ref{dfn11.3} and \ref{dfn11.5} for a precise 
formulation. The prototypical examples are:
\begin{enumerate}
\item $A$ is a finite $C$-algebra (e.g.\ an Azumaya algebra);
\item $A$ is the ring ${\mcal D}(C)$ of differential operators
of $C$, where $C$ is smooth and $\opn{char} \k = 0$; and
\item $A$ is the universal enveloping algebra 
$\mrm{U}(C; L)$ of a finite Lie algebroid $L$ over $C$.
\end{enumerate}
In (1) and (3) there are no regularity assumptions on $C$, 
$A$ or $L$. It is not hard to see that any quotient $A / I$ is also 
a differential $\k$-algebra of finite type. Also if $A_1$ and $A_2$ 
are differential $\k$-algebras of finite type then so is the tensor 
product $A_1 \otimes A_2$. 

The key technical result is the following ``Theorem on the Two 
Filtrations'':

\begin{thm} \label{thm0.5}
Let $A$ be a differential $\k$-algebra of finite type.
Then there is a nonnegative exhaustive filtration
$G = \{ G_{i} A \}_{i \in \mbb{Z}}$ 
of $A$ such that $\opnt{gr}^{G} A$ is a commutative,
finitely generated, connected graded $\k$-algebra.
\end{thm}

Theorem \ref{thm0.5} is proved in Section \ref{sec3}, where it is 
restated as Theorem \ref{thm11.7}. 

Since Van den Bergh's criterion can now be applied, and using 
results from \cite{YZ1}, we obtain:

\begin{cor} 
Let $A$ be a differential $\k$-algebra of finite type. Then $A$ has 
an Auslander rigid dualizing complex $R_A$. For any finite 
$A$-module $M$ one has $\opn{Cdim} M = \opn{GKdim} M$.
\end{cor}

The corollary is proved in Section \ref{sec8}, where it is 
restated as Theorem \ref{thm12.1}. We remind that 
$\opn{GKdim} M$ is the Gelfand-Kirillov dimension of $M$.
The canonical dimension $\opn{Cdim} M$ is defined by
\[ \opn{Cdim} M := 
- \inf\, \{ q \mid \opn{Ext}^q_A(M, R_A) \neq 0 \}
\in {\mbb{Z}} \cup \{ -\infty \}  \]
for a finite $A$-module $M$, and by 
\[ \opn{Cdim} M := 
\opn{sup}\, \{ \opn{Cdim} M' \mid M' \subset M \text{ is finite} \} 
\]
in general. The Auslander property says that 
$\opn{Cdim} M' \leq \opn{Cdim} M$ 
for any finite left or right $A$-module $M$ and any submodule
$M' \subset M$; and it implies that $\opn{Cdim}$ is an exact 
dimension function.

The {\em rigid perverse t-structure} on 
$\msf{D}^{\mrm{b}}_{\mrm{f}}(\cat{Mod} A)$ is defined as follows:
\[ {}^{p}\msf{D}^{\mrm{b}}_{\mrm{f}}(\cat{Mod} A)^{\leq 0} :=
\{ M \in \msf{D}^{\mrm{b}}_{\mrm{f}}(\cat{Mod} A) \mid
\mrm{H}^i \opn{RHom}_A(M, R_A) = 0 
\text{ for all } i < 0 \} \]
and 
\[ {}^{p}\msf{D}^{\mrm{b}}_{\mrm{f}}(\cat{Mod} A)^{\geq 0} :=
\{ M \in \msf{D}^{\mrm{b}}_{\mrm{f}}(\cat{Mod} A) \mid
\mrm{H}^i \opn{RHom}_A(M, R_A) = 0 
\text{ for all } i > 0 \} . \]
The heart 
\[ {}^{p}\msf{D}^{\mrm{b}}_{\mrm{f}}(\cat{Mod} A)^{0} :=
{}^{p}\msf{D}^{\mrm{b}}_{\mrm{f}}(\cat{Mod} A)^{\leq 0} \cap
{}^{p}\msf{D}^{\mrm{b}}_{\mrm{f}}(\cat{Mod} A)^{\geq 0} \]
is called the category of {\em perverse $A$-modules}. It is an 
abelian category, dual to the category 
$\cat{Mod}_{\mrm{f}} A^{\mrm{op}}$ of finite 
$A^{\mrm{op}}$-modules. 

Here is an alternative characterization of the rigid perverse 
t-structure on \linebreak
$\msf{D}^{\mrm{b}}_{\mrm{f}}(\cat{Mod} A)$,
which resembles the original definition in 
\cite{BBD}. For a module $M$ and any integer $i$ define 
$\Gamma_{\msf{M}_i} M$ to be the biggest submodule of 
$M$ with $\opn{Cdim} \leq i$. This is a functor 
$\Gamma_{\msf{M}_i} : \cat{Mod} A \to \cat{Mod} A$, 
and we denote by $\mrm{H}^j_{\msf{M}_i}$ its $j$th right 
derived functor. The next result is 
a special case of Theorem \ref{thm8.6}:

\begin{thm} \label{thm0.4}
Let $A$ be a differential $\k$-algebra of finite type and 
$M \in \msf{D}^{\mrm{b}}_{\mrm{f}}(\cat{Mod} A)$. 
\begin{enumerate}
\item $M \in 
{}^{p}\msf{D}^{\mrm{b}}_{\mrm{f}}(\cat{Mod} A)^{\leq 0}$
iff $\opn{Cdim} \mrm{H}^j M < i$ for all integers $i, j$ such that 
$j > -i$.
\item $M \in 
{}^{p}\msf{D}^{\mrm{b}}_{\mrm{f}}(\cat{Mod} A)^{\geq 0}$
iff $\mrm{H}^j_{\msf{M}_i} M = 0$ for all integers $i, j$ such that 
$j < -i$.
\end{enumerate}
\end{thm}

Interestingly, in a recent preprint \cite{Ka} Kashiwara has 
proved a similar result; see Remark \ref{rem8.2}.

Suppose $C$ is a finitely generated commutative $\k$-algebra. A 
$C$-bimodule $M$ is called a {\em differential $C$-bimodule} if it 
has some bounded below exhaustive filtration 
$F = \{ F_i M \}_{i \in \mbb{Z}}$ by $C$-sub-bimodules, such that 
$\opnt{gr}^{F} M$ is a central $C$-bimodule. This equivalent to 
the condition that the support of the $C^{\mrm{e}}$-module
$M$ is in the diagonal $\Delta(U) \subset U^2$, where
$U := \opn{Spec} C$ (see Proposition \ref{prop4.1}). 

\begin{thm} \label{thm0.3}
Let $C$ be a finitely generated commutative $\k$-algebra and 
$A$ a differential $C$-ring of finite type. Let $R_A$ be the rigid 
dualizing complex of $A$. Then for every $i$ the cohomology 
bimodule $\mrm{H}^i R_A$ is a differential $C$-bimodule.
\end{thm}

This theorem is repeated as Theorem \ref{thm12.9}. One consequence 
is that the rigid dualizing complex $R_A$ localizes on 
$\opn{Spec} C$ (see Corollary \ref{cor8.10}).

Let $A^{\mrm{e}} := A \otimes A^{\mrm{op}}$. It too is a 
differential $\k$-algebra of finite type, so it has a rigid 
dualizing complex $R_{A^{\mrm{e}}}$, and the category 
${}^{p}\msf{D}^{\mrm{b}}_{\mrm{f}}(\cat{Mod} A^{\mrm{e}})^{0}$
of perverse $A^{\mrm{e}}$-modules exists.
By definition the rigid dualizing complex $R_A$ of $A$ is an 
object of 
$\msf{D}^{\mrm{b}}(\cat{Mod} A^{\mrm{e}})$.

\begin{thm} \label{thm0.2}
Let $A$ be a differential $\k$-algebra of finite type. Then the  
rigid dualizing complex $R_A$ is a perverse $A^{\mrm{e}}$-module,
i.e.\ 
$R_A \in 
{}^{p}\msf{D}^{\mrm{b}}_{\mrm{f}}(\cat{Mod} A^{\mrm{e}})^{0}$.
\end{thm}

This theorem is repeated as Theorem \ref{thm12.2} 
in the body of the paper. In addition to being interesting in 
itself, this result is used in \cite{YZ3} to glue rigid dualizing 
complexes on noncommutative ringed schemes -- as perverse 
bimodules.

This paper was previously entitled ``Differential algebras of 
finite type''. 

\medskip \noindent
\textbf{Acknowledgments.}
The authors wish to thank Joseph Bernstein, Sophie 
Chemla, Masaki Kashiwara, Thierry Levasseur, Zinovy Reichstein, 
Paul Smith and Michel Van den Bergh for helpful conversations.

\section{Filtrations of Rings}

By a {\em filtration} of a 
$\k$-algebra $A$ we mean an ascending filtration 
$F = \{ F_{i} A \}_{i \in \mbb{Z}}$ by $\k$-submodules 
such that $1\in F_0 A$ and 
$F_{i} A \cdot F_{j} A \subset F_{i + j} A$.
We shall call $(A, F)$ a filtered $\k$-algebra; but often we shall 
just say that $A$ is a filtered algebra and leave $F$ implicit.

Suppose $(A, F)$ is a filtered $\k$-algebra. 
Given an $A$-module $M$, by an $(A, F)$-filtration of $M$ 
we mean an ascending  filtration 
$F = \{ F_{i} M \}_{i \in \mbb{Z}}$ of $M$
by $\k$-submod\-ules such that 
$F_{i} A \cdot F_{j} M \subset F_{i + j} M$ for all $i$ and $j$.
We call $(M, F)$ a filtered $(A, F)$-module, and allow ourselves 
to drop reference to $F$ when no confusion may arise.

We say the filtration $F$ on $M$ is {\em exhaustive} if 
$M = \bigcup_i F_i M$, 
$F$ is {\em separated} if $0 = \bigcap_i F_i M$,
$F$ is {\em bounded below} if $F_{i_0 - 1} M = 0$ for some integer
$i_0$, and $F$ is {\em nonnegative} if $F_{-1} M = 0$.
The trivial filtration on $M$ is $F_{-1} M := 0$, 
$F_{0} M := M$.

Let us recall some facts about associated graded modules,
and establish some notation. It shall be convenient to use 
the ordered semigroup 
$\mbb{Z} \cup \{ -\infty \}$
where $- \infty < i$ for every $i \in \mbb{Z}$, and 
$i + j := -\infty$ if either $i = -\infty$ or $j = -\infty$. 

Let $(M, F)$ be an exhaustive filtered module. 
The associated graded module is 
\[ \opnt{gr}\, (M, F) = \opnt{gr}^{F} M =  
\boplus_{i \in \mbb{Z}} \opnt{gr}^{F}_i M :=
\boplus_{i \in \mbb{Z}} \frac{F_i M}{F_{i - 1} M} . \]
Given an element $m \in M$ the $F$-degree of $m$ is
\[ \opn{deg}^F(m) := \inf\, \{ i \mid m \in F_i M \} \in \mbb{Z}
\cup \{ -\infty \} . \]
The $F$-symbol of $m$ is 
\[ \opn{symb}^F(m) := m + F_{i - 1} M \in \opnt{gr}^F_i M \]
if $i = \opn{deg}^F(m) \in \mbb{Z}$; and $\opn{symb}^F(m) := 0$
if $\opn{deg}^F(m) = -\infty$.
Thus the homogeneous elements of $\opnt{gr}^F M$ are the 
symbols. 

Recall that the product on the graded algebra 
$\opnt{gr}^F A$ is defined on symbols as follows. Given elements
$a_1, a_2 \in A$ let
$d_i := \opn{deg}^F(a_i)$ and
$\bar{a}_i := \opn{symb}^F(a_i)$. 
If both $d_i > -\infty$ then 
\[ \bar{a}_1 \cdot \bar{a}_2 := 
a_1 \cdot a_2 + F_{d_1 + d_2 - 1} A \in 
\opnt{gr}^F_{d_1 + d_2} A . \]
Otherwise 
$\bar{a}_1 \cdot \bar{a}_2 := 0$.
Similarly one defines a graded $(\opnt{gr}^F A)$-module structure 
on a filtered module $M$.

If $A = \boplus_{i \in \mbb{Z}} A_i$ is a graded algebra then $A$ 
is also filtered, where 
\[ F_i A := \boplus_{j \leq i} A_j . \]
The filtration $F$ is exhaustive and separated. Moreover 
$A \cong \opnt{gr}^{F} A$ as graded algebras. The isomorphism 
sends $a \in A_i$ to its symbol 
$\opnt{symb}^{F}(a) \in \opnt{gr}^{F}_i A$.

\begin{lem} \label{lem10.4}
Suppose the $\k$-algebra $A$ is generated by a sequence of
elements $\{ a_i \}_{i \in I}$, where $I$ is an indexing set
\tup{(}possibly infinite\tup{)}. Given a sequence 
$\{ d_i \}_{i \in I}$ of nonnegative integers,
there is a unique nonnegative exhaustive filtration 
$F = \{ F_d A \}_{d \in \mbb{Z}}$ such that:
\begin{enumerate}
\rmitem{i} For every $d$, $F_d A$ is the $\k$-linear span of the 
products  $a_{j_1} \cdots a_{j_m}$ such that
$d_{j_1} + \cdots + d_{j_m} \leq d$.
\rmitem{ii} The graded algebra $\opnt{gr}^F A$ is generated by 
a sequence of elements $\{ \bar{a}_i \}_{i \in I}$, where 
for every $i \in I$ either $\bar{a}_i = \opnt{symb}^{F}(a_i)$
or $\bar{a}_i = 0$.
\end{enumerate}
\end{lem}

\begin{proof}
Let $\bsym{x} = \{ x_i \}_{i \in I}$ be a sequence of distinct
indeterminates, and let $\k \bra{\bsym{x}}$
be the free associative algebra on these generators. Define
$\phi: \k \bra{\bsym{x}} \to A$
to be the surjection sending $x_i \mapsto a_i$.
Put on $\k \bra{\bsym{x}}$ the grading such that
$\opn{deg}(x_i) = d_i$. 
This induces a filtration 
$F = \{ F_d \k \bra{\bsym{x}} \}_{d \in \mbb{Z}}$
where
\[ F_d \k \bra{\bsym{x}} :=
\boplus_{e \leq d} \k \bra{\bsym{x}}_e . \]
This filtration can now be transferred to $A$ by setting
$F_d A := \phi(F_d \k \bra{\bsym{x}})$. 
Clearly $(A, F)$ is exhaustive and nonnegative, and also 
condition (i) holds. This condition also guarantees uniqueness.

As for condition (ii) consider the surjective graded algebra 
homomorphism
\[ \opnt{gr}^F(\phi): \opnt{gr}^F \k \bra{\bsym{x}} \to
\opnt{gr}^F A . \]
Because of the way the filtration on $\k \bra{\bsym{x}}$ was 
constructed the graded algebra 
$\opnt{gr}^F \k \bra{\bsym{x}}$ is a free algebra on the symbols
$\bar{x}_i := \opn{symb}^F(x_i)$. Define
$\bar{a}_i := \opnt{gr}^F(\phi)(\bar{x}_i)$.
These elements have the required properties.
\end{proof}

Conversely we have the following two lemmas, whose standard
proofs we leave out. 

\begin{lem}
\label{lem10.1}
Let $F=\{ F_d A \}$ be an exhaustive nonnegative filtration 
of $A$, and let $\{ a_i \}_{i \in I}$ be a sequence in $A$.
Denote by $\bar{a}_i  := \opn{symb}^F(a_i)$. Suppose that the 
sequence $\{ \bar{a}_i \}_{i \in I}$ generates $\opnt{gr}^F A$ 
as $\k$-algebra. Then:
\begin{enumerate}
\item $A$ is generated by $\{ a_i \}_{i \in I}$ as $\k$-algebra.
\item Let $d_i := \max\, \{ 0, \opn{deg}^F(a_i) \}$.
Then $F$ coincides with the filtration from Lemma 
\tup{\ref{lem10.4}}. 
\end{enumerate}
\end{lem}

\begin{lem} \label{lem10.3}
Let $(A, F)$ be a nonnegative exhaustive filtered $\k$-algebra and 
let $(M, F)$ be a bounded below exhaustive filtered 
$(A, F)$-module. Suppose 
$\{ a_i \}_{i \in I} \subset A$,  
$\{ b_j \}_{j \in J} \subset A$
and $\{ c_k \}_{k \in K} \subset M$ are sequences satisfying:
\begin{enumerate}
\rmitem{i} The set of symbols 
$\{ \bar{a}_i \}_{i \in I} \cup \{ \bar{b}_j \}_{j \in J}$
generates $\opnt{gr}^{F} A$ as $\k$-algebra.
\rmitem{ii} The set of symbols 
$\{ \bar{c}_k \}_{k \in K}$ generates $\opnt{gr}^{F} M$ as 
$(\opnt{gr}^{F} A)$-module.
\rmitem{iii} For every $i, j$ the symbols
$\bar{a}_i$ and $\bar{b}_j$ commute.
\end{enumerate}
Then for every integer $d$ the $\k$-module $F_d M$ is generated by 
the set of products
\[ \{ a_{i_1} \cdots a_{i_p} b_{j_1} \cdots b_{j_q} c_k \mid
\opn{deg}^F(a_{i_1}) + \cdots 
+ \opn{deg}^F(b_{j_1}) + \cdots 
+ \opn{deg}^F(c_{k}) \leq d \} . \]
\end{lem}

The base field $\k$ is of course trivially filtered.
The filtered $\k$-modules $(M, F)$ form an 
additive  category $\cat{FiltMod} \k$, in which a morphism 
$\phi: (M, F) \to (N, F)$ is a $\k$-linear homomorphism 
$\phi: M \to N$ such that $\phi(F_{i} M) \subset F_{i} N$. 

The {\em Rees module} of $(M, F)$ is 
\[ \opnt{Rees}\, (M, F) = \opnt{Rees}^{F} M := 
\bigoplus_{i \in \mbb{Z}} F_{i} M \cdot t^{i} \subset
M[t] = M \otimes_{\k} \k[t] \]
where $t$ is a central indeterminate of degree $1$. 
We get an additive functor
\[ \opnt{Rees} : \cat{FiltMod} \k \to 
\cat{GrMod} \k[t] \]
where $\cat{GrMod} \k[t]$ is the abelian category of 
graded $\k[t]$-modules and degree $0$ homomorphisms.

For a scalar $\lambda \in \k$ let us denote by
$\opn{sp}_{\lambda}$ the specialization to $\lambda$ of a 
$\k[t]$-module $\wtil{M}$, namely
\[ \opn{sp}_{\lambda} \wtil{M} := 
\wtil{M} / (t - \lambda) \wtil{M} . \]
If $\lambda \neq 0$ then 
\[ \opn{sp}_{\lambda} : \cat{GrMod} \k[t] \to 
\cat{Mod} \k \]
is an exact functor, since $\opn{sp}_{\lambda} \wtil{M}$ 
is isomorphic to the degree zero component of the localization 
$\wtil{M}_{t}$. 
For $\lambda  = 0$ we get a functor 
\[ \opn{sp}_{0}: \cat{GrMod} \k[t] \to 
\cat{GrMod} \k . \]

Given any $(M, F) \in \cat{FiltMod} \k$ one has
\[ \opn{sp}_{0} \opnt{Rees}\, (M, F) \cong \opnt{gr}\, (M, F) = 
\opnt{gr}^{F} M . \]
On the other hand, given a graded $\k[t]$-module $\wtil{M}$ 
there is a filtration $F$ on $M := \opn{sp}_{1} \wtil{M}$ 
defined by
\[ F_{i} M := \opn{Im}
\bigl( \bigoplus\nolimits_{j \leq i} 
\wtil{M}_{j} \to M \bigr) . \]
This is a functor
\[ \opn{sp}_{1}: \cat{GrMod} \k[t] \to 
\cat{FiltMod} \k . \]
If $(M, F)$ is exhaustive then 
\[ \opn{sp}_{1} \opnt{Rees}\, (M, F) \cong  (M, F) .\]
For a graded module $\wtil{M} \in \cat{GrMod} \k[t]$ we have
\[ \opnt{Rees} \opn{sp}_{1} \wtil{M} \cong 
\wtil{M} / \{ t\text{-torsion} \} . \]

If $A$ is a filtered $\k$-algebra then 
$\wtil{A} := \opn{Rees} A$ and
$\bar{A} := \opn{gr} A$ are graded algebras, and we obtain 
corresponding functors $\opnt{Rees}$, $\opn{sp}_{1}$ and 
$\opn{sp}_{0}$ between
$\cat{FiltMod} A$, $\cat{GrMod} \wtil{A}$ and 
$\cat{GrMod} \bar{A}$.

The next lemma says that a filtration can be lifted to the Rees 
ring.

\begin{lem} \label{lem10.7}
Let $F$ be an exhaustive nonnegative filtration of the 
$\k$-algebra $A$, and let
$\wtil{A} := \opnt{Rees}^F A \subset A[t]$. For any integer $i$ 
define
\[ \wtil{F}_i \wtil{A} := \boplus_{j \in \mbb{N}}
(F_{\opn{min} (i, j)} A) \cdot t^j \in A[t] . \]
Then:
\begin{enumerate}
\item $\wtil{A} = \bigcup_i \wtil{F}_i \wtil{A}$,
$\wtil{F}_{-1} \wtil{A} = 0$ and 
$\wtil{F}_i \wtil{A} \cdot \wtil{F}_j \wtil{A}
\subset \wtil{F}_{i + j} \wtil{A}$.
Thus $\wtil{F} = \{ \wtil{F}_i \wtil{A} \}$ is an exhaustive 
nonnegative filtration of the algebra $\wtil{A}$.
\item There is an isomorphism of $\k$-algebras
\[ \opnt{gr}^{\wtil{F}} \wtil{A} \cong
(\opnt{gr}^{F} A) \otimes \k[t]  \]
\tup{(}not respecting degrees\tup{)}. 
\end{enumerate}
\end{lem}

\begin{proof}
(1) Since $F_{-1} A = 0$ we get $\wtil{F}_{-1} \wtil{A} = 0$.
Let $\wtil{A}_j := (F_j A) \cdot t^j$, the $j$th graded component 
of $\wtil{A}$. Then for any $i, j$ there is equality
\[ \wtil{A}_j \cap (\wtil{F}_{i} \wtil{A}) =
(F_{\opn{min} (i, j)} A) \cdot t^j . \]
Hence $\wtil{A}_j \subset \wtil{F}_{j} \wtil{A}$. It remains to 
check the products. For any two pairs of numbers
$(i, k)$ and $(j, l)$ one has
\[ \opn{min} (i, k) + \opn{min} (j, l) \leq
\opn{min} (i + j, k + l) . \]
Therefore
\[ \bigl( (F_{\opn{min} (i, k)} A) \cdot t^k \bigr) \cdot
\bigl( (F_{\opn{min} (j, l)} A) \cdot t^l \bigr) \subset
\bigl( (F_{\opn{min} (i + j, k + l)} A) \cdot t^{k + l} \bigr) . \]
This says that
$\wtil{F}_i \wtil{A} \cdot \wtil{F}_j \wtil{A}
\subset \wtil{F}_{i + j} \wtil{A}$.

\medskip \noindent
(2) We have isomorphisms
\[ \opnt{gr}^{\wtil{F}} \wtil{A} \cong
\boplus_{0 \leq i} \boplus_{i \leq j} 
(\opnt{gr}^{F}_i A) \cdot t^j  \]
and
\[ (\opnt{gr}^{F} A) \otimes \k[t] \cong
\boplus_{0 \leq i} \boplus_{0 \leq j} 
(\opnt{gr}^{F}_i A) \cdot t^j . \]
The isomorphism 
$\opnt{gr}^{\wtil{F}} \wtil{A} \iso (\opnt{gr}^{F} A) \otimes 
\k[t]$
we want is defined on every summand 
$(\opnt{gr}^{F}_i A) \cdot t^j$
by dividing by $t^i$.
\end{proof}

\section{Differential $\k$-Algebras of Finite Type}

Let $C$ be a ring. Recall that a {\em $C$-ring} is a ring $A$ 
together with a ring homomorphism $C \to A$, called the 
structural homomorphism.  
Observe that a $C$-ring is also a $C$-bimodule. 

\begin{dfn} 
\label{dfn11.3}
Suppose $C$ is a commutative $\k$-algebra and $A$ 
is a $C$-ring. A {\em differential $C$-filtration} on
$A$ is a filtration $F = \{ F_{i} A \}_{i \in \mbb{Z}}$ with the 
following properties:
\begin{enumerate}
\rmitem{i} $1 \in F_0 A$ and 
$F_i A \cdot F_j A \subset F_{i + j} A$.
\rmitem{ii} $F_{-1} A = 0$ and $A = \bigcup F_{i} A$.
\rmitem{iii} Each $F_i A$ is a $C$-sub-bimodule. 
\rmitem{iv} The graded ring $\opnt{gr}^{F} A$ 
is a $C$-algebra.
\end{enumerate}
A is called a {\em differential $C$-ring} if it admits 
some differential $C$-filtration.
\end{dfn}

The name ``differential filtration'' signifies the similarity to 
Grothendieck's definition of differential operators; see 
\cite{EGA-IV}.

Note that properties (i) and (iii) imply that the image of the 
structural homomorphism $C \to A$ lies in $F_0 A$, so that 
(iv) makes sense. 

\begin{dfn}
\label{dfn11.5}
Let $C$ be a commutative noetherian $\k$-algebra and $A$ a $C$-ring. 
A {\em differential $C$-filtration of finite type} on $A$ is a 
differential $C$-filtration $F = \{ F_{i} A \}$ such that the 
graded $C$-algebra $\opnt{gr}^{F} A$ is a finite module 
over its center $\mrm{Z}(\opnt{gr}^{F} A)$, and 
$\mrm{Z}(\opnt{gr}^{F} A)$  is a finitely generated $C$-algebra.
We say $A$ is a {\em differential $C$-ring of finite type } if it 
admits some differential $C$-filtration of finite type. 
\end{dfn}

If $C$ is a finitely generated commutative $\k$-algebra and $A$ 
is a differential $C$-ring of finite type then $A$ is also a 
differential $\k$-ring of finite type. In this case we also call 
$A$ a {\em differential $\k$-algebra of finite type}.

\begin{exa} \label{exa11.0}
Let $C$ be a finitely generated commutative $\k$-algebra and $A$ a 
finite $C$-algebra. Then $A$ is a differential $C$-ring of finite 
type. As filtration we can take the trivial filtration 
$F_{-1} A := 0$ and $F_0 A := A$.
\end{exa}

\begin{exa} \label{exa11.1}
Suppose $\opn{char} \k = 0$, $C$ is a smooth commutative 
$\k$-algebra and $A := \mcal{D}(C)$ is the ring of $\k$-linear 
differential operators. Then $A$ is a differential $C$-ring of 
finite type. For filtration we can take the filtration 
$F = \{ F_i A \}$ by order of operator, in which 
$F_{-1} A := 0$, $F_0 A := C$, $F_1 A := C \oplus \mcal{T}(C)$
and $F_{i + 1} A := F_i A \cdot F_1 A$ for $i \geq 1$. Here
$\mcal{T}(C) := \opn{Der}_{\k}(C)$, the module of derivations.
\end{exa}

\begin{exa} \label{exa11.2}
A special case of Example \ref{exa11.1} is when 
$C := \k[x_1, \ldots, x_n]$, a polynomial ring. Then $A$ is called 
the $n$th Weyl algebra. Writing 
$y_i := \frac{\partial}{\partial x_i}$
the algebra $A$ is generated by the $2n$ elements 
$x_1, \ldots, x_n, y_1, \ldots, y_n$, with relations
$[x_i, x_j] = [y_i, y_j] = [y_i, x_j] = 0$ 
for $i \neq j$ and $[y_i, x_i] = 1$.
In addition to the filtration $F$ above there is also a 
differential $\k$-filtration of finite type 
$G = \{ G_i A \}$ where $G_{-1} A := 0$, $G_0 A := \k$, 
$G_1 A := \k + (\sum_i \k \cdot x_i) + (\sum_i \k \cdot y_i)$
and $G_{i + 1} A := G_i A \cdot G_1 A$ for $i \geq 1$.
\end{exa}

\begin{exa} \label{exa11.4}
This example generalizes Example \ref{exa11.2}.
Let $C$ be a finitely generated commutative $\k$-algebra (not 
necessarily smooth, and $\opn{char} \k$ arbitrary), and let $L$ be 
a finite $C$-module (not necessarily projective). Suppose $L$ has 
a $\k$-linear Lie bracket $[-, -]$.
The module of derivations $\mcal{T}(C) := \opn{Der}_{\k}(C)$
is also a $\k$-Lie algebra. Suppose 
$\alpha : L \to \mcal{T}(C)$
is a $C$-linear Lie homomorphism, namely 
$\alpha(c \xi) = c \alpha(\xi)$ and 
$\alpha([\xi, \zeta]) = [\alpha(\xi), \alpha(\zeta)]$
for all $c \in C$ and $\xi, \zeta \in L$. $L$ is then called a 
{\em Lie algebroid} or a {\em Lie-Rinehart algebra}
(cf.\ \cite{Ch1} or \cite{Ri}). The ring of 
generalized differential operators $\mcal{D}(C; L)$, also called 
the universal enveloping algebra and denoted $\mrm{U}(C; L)$, is 
defined as follows. Choose $\k$-algebra generators 
$c_1, \ldots, c_p$ for $C$ and $C$-module generators 
$l_1, \ldots, l_q$ for $L$. Let
\[ \k \bra{\bsym{x}, \bsym{y}} :=
\k \bra{x_{1}, \ldots, x_{p}, y_{1}, \ldots, y_{q}} \] 
be the free associative algebra.
We have a ring surjection 
$\phi_0 : \k \bra{\bsym{x}} \to C$ with $\phi_0(x_i) := c_i$. 
Let $I_0 := \opn{Ker}(\phi_0)$. Next there is a surjection of 
$\k \bra{\bsym{x}}$-modules
\[ \phi_1 : \k \bra{\bsym{x}}^q = 
\boplus_{j=1}^q \k \bra{\bsym{x}} \cdot y_j \to L \]
with $\phi_1(y_j) := l_j$. Define 
$I_1 := \opn{Ker}(\phi_1) \subset \k \bra{\bsym{x}}^q$. 
For any $i, j$ choose polynomials $f_{i, j}(\bsym{x})$ and
$g_{i, j, k}(\bsym{x})$ such that
$[l_i, l_j] = \sum_k g_{i, j, k}(\bsym{c}) l_k \in L$
and
$\alpha(l_i)(c_j) = f_{i, j}(\bsym{c}) \in C$.
Now define
\[ \mrm{U}(C; L) := \frac{\k \bra{\bsym{x}, \bsym{y}}}
{I} \]
where $I \subset \k \bra{\bsym{x}, \bsym{y}}$ 
is the two-sided ideal generated by $I_0$, $I_1$ and the 
polynomials 
$[y_i, y_j] - \sum_k g_{i, j, k}(\bsym{x}) y_k$
and
$[y_i, x_j] - f_{i, j}(\bsym{x})$.

The ring $\mrm{U}(C; L)$ has the following universal property: 
given any ring $D$, any ring homomorphism $\eta_0 : C \to D$ and 
any $C$-linear Lie homomorphism $\eta_1 : L \to D$ satisfying
$[\eta_1(l), \eta_0(c)] = \eta_0(\alpha(l)(c))$,
there is a unique ring homomorphism $\eta : \mrm{U}(C; L) \to D$
through which $\eta_0$ and $\eta_1$ factor.

Put on $\k \bra{\bsym{x}, \bsym{y}}$ the filtration $F$ such 
that $\opn{deg}^F(x_i) = 0$ and $\opn{deg}^F(y_j) = 1$. Let $F$ be 
the filtration induced on $\mrm{U}(C; L)$ by the surjection
$\phi : \k \bra{\bsym{x}, \bsym{y}} \surj
\mrm{U}(C; L)$.
Then 
$\opnt{gr}^{F} \mrm{U}(C; L)$ is a commutative $C$-algebra, 
generated by the elements 
$\bar{l}_j := \opnt{gr}^{F}(\phi)(y_j)$, $j \in \{ 1, \ldots, q \}$.
We see that $\mrm{U}(C; L)$ is a differential $C$-ring of finite 
type. If 
$C = \k[x_1, \ldots, x_n]$ and $L = \mcal{T}(C)$ then we are in the 
situation of Example \ref{exa11.2}. If $C = \k$ then 
$\mrm{U}(C; L) = \mrm{U}(L)$ is the usual universal enveloping 
algebra of the Lie algebra $L$.
\end{exa}

\begin{lem}[{\cite[Theorem 8.2]{ATV}}]
\label{lem11.4}
Suppose $A = \boplus_{i \in \mbb{N}} A_i$ 
is a graded $\k$-algebra and 
$t \in A$ is a central homogeneous element of positive degree. The 
following are equivalent:
\begin{enumerate}
\rmitem{i} $A$ is left noetherian. 
\rmitem{ii} $A / (t)$ is left noetherian. 
\end{enumerate}
\end{lem}

The next proposition follows almost directly from the definition 
and Lemma \ref{lem11.4}. 

\begin{prop} \label{prop11.3}
If $A$ is a differential $\k$-algebra
of finite type then it is a noetherian finitely generated 
$\k$-algebra.
\end{prop}

The class of differential $\k$-algebras
of finite type is closed under tensor products as we now show.

\begin{prop} \label{prop11.6}
Let $C_1$ and $C_2$ be noetherian commutative $\k$-algebras, and 
assume that $C_{1} \otimes C_{2}$ is also noetherian. 
Let $A_i$ be a differential $C_{i}$-ring of finite type for $i=1,2$. 
Then $A_{1} \otimes A_{2}$ is a differential 
$(C_{1} \otimes C_{2})$-ring of finite type. 
\end{prop}

\begin{proof}
Choose differential filtrations of finite type
$\{ F_{n} A_{1} \}$ and $\{ F_{n} A_{2} \}$
of $A_{1}$ and $A_{2}$. Define a filtration on 
$A_{1} \otimes A_{2}$ as follows:
\[ F_{n} (A_{1} \otimes A_{2}) := \sum_{l + m = n}
F_{l} A_{1} \otimes F_{m} A_{2} . \]
Then 
\[ \opnt{gr}^{F} (A_{1} \otimes A_{2}) \cong
(\opnt{gr}^{F} A_{1}) \otimes (\opnt{gr}^{F} A_{2}) \]
as graded rings. Since $\opnt{gr}^F A_1$ and $\opnt{gr}^F A_2$
are finite modules over their centers it follows that 
$(\opnt{gr}^{F} A_{1}) \otimes (\opnt{gr}^{F} A_{2})$
is a finite module over its center.
\end{proof}

\section{The Theorem on the Two Filtrations}
\label{sec3}

The next theorem generalizes the case of the $n$th 
Weyl algebra and its two filtrations (see Examples \ref{exa11.1} 
and \ref{exa11.2} above). McConnell and Stafford also considered 
such filtrations, and our result extends their 
\cite[Corollary 1.7]{MS}. 
The basic idea is attributed in \cite{MS} to Bernstein. 

We recall that a graded $\k$-algebra $A$ is called 
{\em connected} if $A = \boplus_{i \in \mbb{N}} A_i$, 
$A_0 = \k$ and each $A_i$ is a finite $\k$-module.

\begin{thm} \label{thm11.7}
Let $A$ be a $\k$-algebra. Assume $A$ has a differential 
$\k$-filtration of finite type $F = \{ F_{i} A \}_{i \in \mbb{Z}}$.
Then there is a nonnegative exhaustive $\k$-filtration
$G = \{ G_{i} A \}_{i \in \mbb{Z}}$ 
such that $\opnt{gr}^{G} A$ is a commutative,
finitely generated, connected graded $\k$-algebra.
\end{thm}

Observe that $G$ is also a differential filtration of finite 
type on $A$. 
As mentioned in the introduction Theorem \ref{thm11.7} will be 
used to prove the existence of an Auslander rigid dualizing 
complex over $A$.

The following easy lemma will be used often in the proof of 
the theorem.

\begin{lem} \label{lem11.5}
Let $F = \{ F_i A \}$ be a nonnegative exhaustive filtration of 
$A$ and let $a_1, a_2 \in A$ be two elements. Define
$\bar{a}_i := \opn{symb}^F(a_i) \in \opnt{gr}^{F} A$
and
$d_i := \opn{deg}^F(a_i) \in \mbb{N} \cup \{ -\infty \}$.
Then the commutator $[\bar{a}_1, \bar{a}_2] = 0$ if and only if
\[ \opn{deg}^F([a_1, a_2]) \leq d_1 + d_2 - 1 . \]
\end{lem}

\begin{proof}[Proof of Theorem \tup{\ref{thm11.7}}]
Step 1. 
Write $\bar{A} := \opnt{gr}^{F} A$. Then the center 
$\mrm{Z}(\bar{A})$ is a graded, finitely generated, commutative 
$\k$-algebra, and $\bar{A}$ is a finite $\mrm{Z}(\bar{A})$-module.
For any element $a \in A$ we write 
$\bar{a} := \opn{symb}^F(a) \in \opnt{gr}^F A$. 

Let $d_1 \in \mbb{N}$ be large enough such that 
$\mrm{Z}(\bar{A})$ is generated 
as $\mrm{Z}(\bar{A})_{0}$-algebra by finitely many elements of 
degrees $\leq d_1$, and $\bar{A}$ is generated as 
$\mrm{Z}(\bar{A})$-module by finitely many elements of 
degrees $\leq d_1$. 

Choose nonzero elements 
$a_{1}, \ldots, a_{m} \in F_{0} A \cong
\bar{A}_{0}$ such that their symbols 
$\bar{a}_{1}, \ldots, \bar{a}_{m}$ are in 
$\mrm{Z}(\bar{A})_{0}$, and they generate $\mrm{Z}(\bar{A})_{0}$ as 
$\k$-algebra. Next choose elements 
$b_{1}, \ldots, b_{n} \in F_{d_1} A - F_0 A$ 
such that the symbols 
$\bar{b}_{1}, \ldots, \bar{b}_{n}$ are in 
$\mrm{Z}(\bar{A})$, and they generate $\mrm{Z}(\bar{A})$ as 
$\mrm{Z}(\bar{A})_{0}$-algebra. 
Finally choose nonzero elements 
$c_{1}, \ldots, c_{p} \in F_{d_1} A$ 
such that the symbols 
$\bar{c}_{1}, \ldots, \bar{c}_{p}$ generate 
$\boplus_{j=0}^{d_1} \bar{A}_j$ as $\mrm{Z}(\bar{A})_0$-module. 
This implies that 
$\bar{c}_{1}, \ldots, \bar{c}_{p}$ generate $\bar{A}$ as 
a $\mrm{Z}(\bar{A})$-module.

The symbols 
$\bar{a}_{1}, \ldots, \bar{a}_{m},
\bar{b}_{1}, \ldots, \bar{b}_{n},
\bar{c}_{1}, \ldots, \bar{c}_{p}$
generate $\bar{A}$ as $\k$-algebra, so by \linebreak
Lemma \ref{lem10.1} the elements
$a_{1}, \ldots, c_{m}, b_{1}, \ldots, b_{n},
c_{1}, \ldots, c_{p}$
generate $A$ as $\k$-algebra. Let 
\[ \k \bra{\bsym{x}, \bsym{y}, \bsym{z}} :=
\bra{x_{1}, \ldots, x_{m}, y_{1}, \ldots, y_{n}, 
z_{1}, \ldots, z_{p}} \] 
be the free associative algebra, and define a surjective ring 
homomorphism
$\phi: \k \bra{\bsym{x}, \bsym{y}, \bsym{z}} \to A$
by sending $x_i \mapsto a_i$, $y_i \mapsto b_i$
and $z_i \mapsto c_i$. We are now in the situation of Lemma 
\ref{lem10.1}. The free algebra 
$\k \bra{\bsym{x}, \bsym{y}, \bsym{z}}$
also has a filtration $F$, where
$\opn{deg}^F(x_i) := \opn{deg}^F(a_i)$, 
$\opn{deg}^F(y_i) := \opn{deg}^F(b_i)$
etc., and $\phi$ is a strict surjection, meaning that 
$F_i(A) = \phi(F_i  \k \bra{\bsym{x}, \bsym{y}, \bsym{z}})$.

Let us denote substitution by 
$f(\bsym{a}, \bsym{b}, \bsym{c}) := 
\phi(f(\bsym{x}, \bsym{y}, \bsym{z}))$. 
Consider the subrings 
\[ \k \bra{\bsym{a}} \subset \k \bra{\bsym{a}, \bsym{b}} \subset A 
= \k \bra{\bsym{a}, \bsym{b}, \bsym{c}} \]
with filtrations $F$ induced by the inclusions into $A$.
Warning: these filtrations might differ from the 
filtrations induced by 
$\phi: \k \bra{\bsym{x}} \surj \k \bra{\bsym{a}}$ and 
$\phi: \k \bra{\bsym{x}, \bsym{y}} 
\surj \k \bra{\bsym{a}, \bsym{b}}$
respectively. 

We observe that the commutators $[a_i, a_j] = 0$ for all $i, j$,
since $a_i = \bar{a}_i \in Z(\bar{A})_0$. This also says that
$[\bar{a}_i, \bar{b}_j] = 0$, so according to Lemma \ref{lem11.5} 
we get
\[ [a_{i}, b_{j}] \in F_{\opn{deg}^F(b_{j}) - 1} A \]
for all $i,j$. Therefore by Lemma \ref{lem10.3}, applied to the 
filtered $\k$-algebra $\k \bra{\bsym{a}}$ and the filtered
$\k \bra{\bsym{a}}$-module
$F_{d_1} A$, we see that there are noncommutative
polynomials 
$f^1_{i, j, k}(\bsym{x}) \in \k \bra{\bsym{x}}$ 
such that 
\begin{equation} 
\label{eqn11.8}
\begin{gathered}
\opn{deg}^F(f^1_{i, j, k}(\bsym{x})) + \opn{deg}^F(c_k) \leq
\opn{deg}^F(b_j) - 1 \\
\text{and} \\
[a_{i}, b_{j}] = \sum_{k} f^1_{i, j, k}(\bsym{a}) \cdot c_{k} . 
\end{gathered}
\end{equation}
Note that either $f^1_{i, j, k}(\bsym{x}) \neq 0$, in which case 
$\opn{deg}^F(f^1_{i, j, k}(\bsym{x})) = 0$; 
or $f^1_{i, j, k}(\bsym{x}) = 0$ and then 
$\opn{deg}^F(f^1_{i, j, k}(\bsym{x})) = -\infty$. 
The choice $f^1_{i, j, k}(\bsym{x}) = 0$ is of course required when 
$\opn{deg}^F(c_k) \geq \opn{deg}^F(b_j)$.

Likewise 
$[b_{i}, b_{j}] \in F_{\opn{deg}^F(b_{i}) + 
\opn{deg}^F(b_{j}) - 1} A$,
so by Lemma \ref{lem10.3}, applied to the filtered $\k$-algebra
$\k \bra{\bsym{a}, \bsym{b}}$ and the filtered
$\k \bra{\bsym{a}, \bsym{b}}$-module
$A$, we see that there are noncommutative polynomials 
$f^2_{i, j, k}(\bsym{x})$ and
$g^2_{i, j, k}(\bsym{y})$
such that 
\begin{equation} \label{eqn11.9}
\begin{gathered}
\opn{deg}^F(f^2_{i, j, k}(\bsym{x})) + 
\opn{deg}^F(g^2_{i, j, k}(\bsym{y})) +
\opn{deg}^F(c_k) \leq
\opn{deg}^F(b_i) + \opn{deg}^F(b_j) - 1 \\
\text{and} \\
[b_{i}, b_{j}] = 
\sum_{k} f^2_{i, j, k}(\bsym{a}) \cdot 
g^2_{i, j, k}(\bsym{b}) \cdot c_{k} .
\end{gathered}
\end{equation}

Similarly there are polynomials 
$f^3_{i, j, k}(\bsym{x})$ such that
\[ \begin{gathered}
\opn{deg}^F(f^3_{i, j, k}(\bsym{x})) + \opn{deg}^F(c_k) \leq
\opn{deg}^F(c_j) - 1 \\
\text{and} \\
[a_{i}, c_{j}] = 
\sum_{k} f^3_{i, j, k}(\bsym{a}) \cdot c_{k} ,
\end{gathered} \]
and there are polynomials 
$f^4_{i, j, k}(\bsym{x})$
and
$g^4_{i, j, k}(\bsym{y})$
such that 
\[ \begin{gathered}
\opn{deg}^F(f^4_{i, j, k}(\bsym{x})) + 
\opn{deg}^F(g^4_{i, j, k}(\bsym{y})) +
\opn{deg}^F(c_k) \leq
\opn{deg}^F(b_i) + \opn{deg}^F(c_j) - 1 \\
\text{and} \\
[b_{i}, c_{j}] =
\sum_{k} f^4_{i, j, k}(\bsym{a}) \cdot
g^4_{i, j, k}(\bsym{b}) \cdot c_{k} . 
\end{gathered} \]

The same idea applies to $c_i c_j$:
there are polynomials $f^5_{i, j, k}(\bsym{x})$ 
and
$g^5_{i, j, k}(\bsym{y})$
such that 
\[ \begin{gathered}
\opn{deg}^F(f^5_{i, j, k}(\bsym{x})) + 
\opn{deg}^F(g^5_{i, j, k}(\bsym{y})) +
\opn{deg}^F(c_k) \leq
\opn{deg}^F(c_i) + \opn{deg}^F(c_j) \\
\text{and} \\
c_{i} \cdot c_{j} =
\sum_{k} f^5_{i, j, k}(\bsym{a}) \cdot
g^5_{i, j, k}(\bsym{b}) \cdot c_{k} . 
\end{gathered} \]

Let $G$ be the standard grading on $\k \bra{\bsym{x}}$, namely
$\opn{deg}^G(x_i) := 1$. This induces a filtration $G$. 
Define 
\[ e_0 :=  \opn{max}\, \{ 0, \opn{deg}^G(f^l_{i, j, k}(\bsym{x})) \}
, \] 
$e_1 := e_0 + 1$ and $e_2 := e_0 + e_1 + 1$. 

Put on the free algebra $\k \bra{\bsym{x}, \bsym{y}, \bsym{z}}$
a new grading $G$ by declaring
\[ \begin{aligned}
\opn{deg}^G (y_{i}) & := e_2 \opn{deg}^F(b_{i}) \\
\opn{deg}^G (z_{i}) & := e_2 \opn{deg}^F(c_{i}) + e_1 ,
\end{aligned} \]
and keeping $\opn{deg}^G (x_{i}) = 1$ as above. 
We get a new filtration $G$ on
$\k \bra{\bsym{x}, \bsym{y}, \bsym{z}}$.
Using this we obtain a new filtration $G$ on $A$ with
\[ G_{i} A := 
\phi(G_{i} \k \bra{\bsym{x}, \bsym{y}, \bsym{z}}) . \]

\medskip \noindent Step 2.
Now we verify that the filtration $G$ has the required properties. 
Since the filtration $G$ on 
$\k \bra{\bsym{x}, \bsym{y}, \bsym{z}}$
is nonnegative exhaustive, and $\phi$ is a strict surjection, 
it follows that the filtration $G$ on $A$ is also 
nonnegative exhaustive. The rest requires some work, 
and in order to simplify our notation we 
are going to ``recycle'' the expressions 
$\bar{A}$, $\bar{a}_i$ etc. From here on we define
$\bar{A} := \opnt{gr}^{G} A$.
We have a surjective graded $\k$-algebra homomorphism
\[ \bar{\phi} := \opnt{gr}^{G}(\phi) :
\opnt{gr}^{G} \k \bra{\bsym{x}, \bsym{y}, \bsym{z}} \to
\opnt{gr}^{G} A = \bar{A} . \]
Let $\bar{x}_i := \opn{symb}^G(x_i)$,
$\bar{y}_i := \opn{symb}^G(y_i)$, etc. Then
\[ \opnt{gr}^{G} \k \bra{\bsym{x}, \bsym{y}, \bsym{z}} =
\k \bra{\bar{\bsym{x}}, \bar{\bsym{y}}, \bar{\bsym{z}}} :=
\k \bra{\bar{x}_1, \ldots, \bar{x}_m, \bar{y}_1, \ldots, 
\bar{y}_n, \bar{z}_1, \ldots, \bar{z}_p} \]
which is also a free algebra. Define
$\bar{a}_{i} := \bar{\phi}(\bar{x}_{i})$, 
$\bar{b}_{i} := \bar{\phi}(\bar{y}_{i})$ and
$\bar{c}_{i} := \bar{\phi}(\bar{z}_{i})$.
Observe that either 
$\opn{deg}^G(a_{i}) = \opn{deg}^G(x_{i})$, in which case 
$\bar{a}_{i} = \opn{symb}^G(a_i)$, and it is a nonzero 
element of $\bar{A}_{\opn{deg}^G(a_{i})}$; or 
$\opn{deg}^G(a_{i}) < \opn{deg}^G(x_{i})$, and then 
$\bar{a}_{i} = 0$. Similar statements hold for $\bar{b}_{i}$ and 
$\bar{c}_{i}$.

Since $\bar{\phi}$ is surjective we see that 
$\bar{A}$ is generated as $\k$-algebra by the elements
$\bar{a}_{1}, \ldots, \bar{a}_{m}, \bar{b}_{1}, \ldots, 
\bar{b}_{n}, \bar{c}_{1},$ 
$\ldots, \bar{c}_{p}$.
These elements are either of positive degree or are $0$, and hence 
$\bar{A}$ is connected graded. 
We claim that $\bar{A}$ is commutative.

We know already that $[\bar{a}_{i}, \bar{a}_{j}] = 0$.
Let us check that $[\bar{a}_{i}, \bar{b}_{j}] = 0$. 
If either 
$\opn{deg}^G(a_i) < \opn{deg}^G(x_i)$ or 
$\opn{deg}^G(b_j) < \opn{deg}^G(y_j)$ then 
$\bar{a}_i \bar{b}_j = \bar{b}_j \bar{a}_i = 0$. 
Otherwise 
$\bar{a}_i = \opn{symb}^G(a_i)$ and 
$\bar{b}_j = \opn{symb}^G(b_j)$. 
By formula (\ref{eqn11.8}) we have
\[ \begin{aligned}
\opn{deg}^G([a_{i}, b_{j}]) & \leq
\opn{max}\, 
\{ \opn{deg}^G(f^1_{i, j, k}(\bsym{a})) + \opn{deg}^G(c_{k}) \} \\
& \leq \opn{max}\, 
\{ \opn{deg}^G(f^1_{i, j, k}(\bsym{x})) + \opn{deg}^G(c_{k}) \}
. 
\end{aligned} \]
For any $k$ such that $f^1_{i, j, k}(\bsym{x}) \neq 0$ we have
$\opn{deg}^F(c_{k}) \leq \opn{deg}^F(b_{j}) - 1$, and then
\[ \opn{deg}^G(c_{k}) \leq \opn{deg}^G(z_{k}) = 
e_2 \opn{deg}^{F}(c_{k}) + e_1 
\leq e_2 (\opn{deg}^{F}(b_{j}) - 1) + e_1 . \]
Also
$\opn{deg}^G(f^1_{i, j, k}(\bsym{x})) \leq e_0.$
Because 
\[ \opn{deg}^G(b_{j}) = \opn{deg}^G(y_{j}) = 
e_2 \opn{deg}^{F}(b_{j}) \]
and 
\[ \opn{deg}^G(a_{i}) = \opn{deg}^G(x_{i}) = 1 , \] 
we get 
\[ \begin{aligned}
\opn{deg}^G([a_{i}, b_{j}]) & \leq 
e_0 + (e_2 (\opn{deg}^{F}(b_{j}) - 1) + e_1) \\
& = \opn{deg}^G(a_{i}) + \opn{deg}^G(b_{j}) - 2 .
\end{aligned} \]
Using Lemma \ref{lem11.5} we conclude that
$[\bar{a}_{i}, \bar{b}_{j}] = 0$.

Next let's consider the commutator $[\bar{b}_{i}, \bar{b}_{j}]$.
If either 
$\opn{deg}^G(b_i) < \opn{deg}^G(y_i)$ or 
$\opn{deg}^G(b_j) < \opn{deg}^G(y_j)$ then 
$\bar{b}_i \bar{b}_j = \bar{b}_j \bar{b}_i = 0$. 
Otherwise 
$\bar{b}_i = \opn{symb}^G(b_i)$ and 
$\bar{b}_j = \opn{symb}^G(b_j)$. 
By formula (\ref{eqn11.9}), if 
$f^2_{i, j, k}(\bsym{x}) \neq 0$
then 
\[ \opn{deg}^F(g^2_{i, j, k}(\bsym{y})) + 
\opn{deg}^F(c_{k}) \leq \opn{deg}^F(b_i) + \opn{deg}^F(b_j) - 1 . 
\]
Also
\[ \begin{aligned}
\opn{deg}^G(g^2_{i, j, k}(\bsym{y})) & = 
e_2 \opn{deg}^F(g^2_{i, j, k}(\bsym{y})) \\
\opn{deg}^G(c_k) & \leq \opn{deg}^G(z_k) =
e_2 \opn{deg}^F(c_k) + e_1 . 
\end{aligned} \]
Therefore, looking only at indices $k$ s.t.\ 
$f^2_{i, j, k}(\bsym{x}) \neq 0$, we obtain
\[ \begin{aligned}
\opn{deg}^G([b_{i}, b_{j}]) & \leq
\opn{max}\, 
\{ \opn{deg}^G(f^2_{i, j, k}(\bsym{a})) + 
\opn{deg}^G(g^2_{i, j, k}(\bsym{b})) + 
\opn{deg}^G(c_{k}) \} \\
& \leq \opn{max}\, 
\{ \opn{deg}^G(f^2_{i, j, k}(\bsym{x})) + 
\opn{deg}^G(g^2_{i, j, k}(\bsym{y})) + 
\opn{deg}^G(c_{k}) \} \\
& \leq e_0 + \opn{max}\, 
\{ e_2 \opn{deg}^F(g^2_{i, j, k}(\bsym{y}))
+ e_2 \opn{deg}^F(c_k) + e_1 \} \\
& = e_2 \opn{max}\, 
\{ \opn{deg}^F(g^2_{i, j, k}(\bsym{y}))
+ \opn{deg}^F(c_k) + 1 \} + (e_0 + e_1 - e_2) \\
& \leq e_2 (\opn{deg}^F(b_i) + \opn{deg}^F(b_j)) - 1 \\
& = \opn{deg}^G(b_i) + \opn{deg}^G(b_j) - 1 .  
\end{aligned} \]
So according to Lemma \ref{lem11.5} we conclude that
$[\bar{b}_{i}, \bar{b}_{j}] = 0$.

The calculation for the other commutators is similar.

Finally we show that, amusingly, $\bar{c}_{i} \bar{c}_{j} = 0$.
If either 
$\opn{deg}^G(c_i) < \opn{deg}^G(z_i)$ or 
$\opn{deg}^G(c_j) < \opn{deg}^G(z_j)$ then automatically
$\bar{c}_i \bar{c}_j = 0$. 
Otherwise 
$\bar{c}_i = \opn{symb}^G(c_i)$ and 
$\bar{c}_j = \opn{symb}^G(c_j)$. 
For any $k$ such that $f^5_{i, j, k}(\bsym{x}) \neq 0$
one has 
\[ \opn{deg}^F(g^5_{i, j, k}(\bsym{y})) + 
\opn{deg}^F(c_{k}) \leq \opn{deg}^F(c_i) + \opn{deg}^F(c_j) . 
\]
Therefore, looking only at indices $k$ s.t.\ 
$f^5_{i, j, k}(\bsym{x}) \neq 0$, we obtain
\[ \begin{aligned}
\opn{deg}^G(c_{i} c_{j}) & \leq
\opn{max}\, 
\{ \opn{deg}^G(f^5_{i, j, k}(\bsym{a})) + 
\opn{deg}^G(g^5_{i, j, k}(\bsym{b})) + 
\opn{deg}^G(c_{k}) \} \\
& \leq \opn{max}\, 
\{ \opn{deg}^G(f^5_{i, j, k}(\bsym{x})) + 
\opn{deg}^G(g^5_{i, j, k}(\bsym{y})) + 
\opn{deg}^G(c_{k}) \} \\
& \leq e_0 + \opn{max}\, 
\{ e_2 \opn{deg}^F(g^5_{i, j, k}(\bsym{y}))
+ e_2 \opn{deg}^F(c_k) + e_1 \} \\
& = e_2 \opn{max}\, 
\{ \opn{deg}^F(g^5_{i, j, k}(\bsym{y}))
+ \opn{deg}^F(c_k) \} + (e_0 + e_1) \\
& \leq e_2 (\opn{deg}^F(c_i) + \opn{deg}^F(c_j)) + (e_0 + e_1) \\
& = \opn{deg}^G(c_i) + \opn{deg}^G(c_j) - 1 .  
\end{aligned} \]
So by definition of the product in $\bar{A}$ we get
$\bar{c}_{i} \bar{c}_{j} = 0$.
\end{proof}

\begin{prop} \label{prop11.2}
In the situation of Theorem \tup{\ref{thm11.7}} 
assume the $\k$-algebra 
$A$ is graded, and also every $\k$-submodule $F_i A$ is graded. 
Then the filtration $G$ can be chosen such that 
every $\k$-submodule $G_i A$ is graded. 
\end{prop}

\begin{proof}
Simply choose the generators 
$a_1, \ldots, b_1, \ldots, c_1, \ldots, c_p \in A$
used in the proof to be homogeneous.
\end{proof}

\section{Review of Dualizing Complexes}
\label{sec5}

For a $\k$-algebra $A$ we denote by $A^{\mrm{op}}$ 
the opposite algebra, and by 
$A^{\mrm{e}} := A \otimes A^{\mrm{op}}$ 
the enveloping algebra. Recall that an $A$-module 
means a left $A$-module. With this convention a right $A$-module 
is an $A^{\mrm{op}}$-module, and an $A$-bimodule is an 
$A^{\mrm{e}}$-module. 

In this section we review the definition of dualizing complexes over 
rings and related concepts. 

Let $\cat{Mod} A$ be the category of $A$-modules, and
let $\cat{Mod}_{\mrm{f}} A$ be the full subcategory of finite 
(i.e.\ finitely generated) modules. The latter is abelian when $A$ 
is left noetherian. Let $\msf{D}(\cat{Mod} A)$ be the derived 
category of $A$-modules. The full subcategory of bounded complexes 
is denoted by $\msf{D}^{\mrm{b}}(\cat{Mod} A)$, the full 
subcategory of complexes with finite cohomologies is denoted by 
$\msf{D}_{\mrm{f}}(\cat{Mod} A)$, and their intersection is 
$\msf{D}^{\mrm{b}}_{\mrm{f}}(\cat{Mod} A)$.

\begin{dfn}[\cite{Ye1}, \cite{YZ1}] \label{dfn2.1}
Let $A$ be a left noetherian $\k$-algebra and $B$ a right noetherian
$\k$-algebra. A complex 
$R \in \msf{D}^{\mrm{b}}(\cat{Mod}\, (A \otimes B^{\mrm{op}}))$ 
is called a {\em dualizing complex over $(A,B)$} if 
it satisfies the following three conditions:
\begin{enumerate}
\rmitem{i} $R$ has finite injective dimension over $A$ and over 
$B^{\mrm{op}}$.
\rmitem{ii} $R$ has finite cohomology modules over $A$ 
and over $B^{\mrm{op}}$.
\rmitem{iii} The canonical morphisms 
$B \to \opn{RHom}_A(R,R)$ in 
$\msf{D}(\cat{Mod} B^{\mrm{e}})$, and 
$A \to \opn{RHom}_{B^{\mrm{op}}}(R,R)$ 
in $\msf{D}(\cat{Mod} A^{\mrm{e}})$, are both isomorphisms.
\end{enumerate}
In the case $A = B$ we say $R$ is a {\em dualizing complex over} 
$A$.
\end{dfn}

Whenever we refer to a dualizing complex over $(A, B)$ in the 
paper we tacitly assume that $A$ is left noetherian and $B$ is 
right noetherian.

\begin{rem}
There are many non-isomorphic dualizing complexes over a given 
$\k$-algebra $A$. The isomorphism classes of dualizing complexes 
are parameterized by the 
derived Picard group $\opn{DPic}(A)$, whose elements are the
isomorphism classes of two-sided tilting complexes. See 
\cite{Ye4} and \cite{MY}.
\end{rem}

Here are two easy examples.

\begin{exa}
Suppose $A$ is a Gorenstein noetherian ring, namely the bimodule 
$R := A$ has finite injective dimension as left and right module. 
Then $R$ is a dualizing complex over $A$.
\end{exa}

\begin{exa}
If $A$ is a finite $\k$-algebra then the bimodule
$A^*:= \opn{Hom}_{\k}(A, \k)$
is a dualizing complex over $A$. In fact it is a rigid dualizing 
complex (see Definition \ref{dfn2.3}). 
\end{exa}

\begin{dfn} \label{dfn2.7}
Let $R$ be a dualizing complex over $(A, B)$. 
The {\em duality functors induced by $R$}  
are the contravariant functors
\[ \mrm{D} := \opn{RHom}_A(-, R): \msf{D}(\cat{Mod} A)
\to \msf{D}(\cat{Mod} B^{\mrm{op}}) \]
and 
\[ \mrm{D}^{\mrm{op}} := \opn{RHom}_{B^{\mrm{op}}}(-, R): 
\msf{D}(\cat{Mod} B^{\mrm{op}}) \to \msf{D}(\cat{Mod} A). \]
\end{dfn}

By \cite[Proposition 1.3]{YZ1} the functors $\mrm{D}$ and 
$\mrm{D}^{\mrm{op}}$
are a duality (i.e.\ an anti-equivalence) of triangulated 
categories between $\msf{D}_{\mrm{f}}(\cat{Mod} A)$ and 
$\msf{D}_{\mrm{f}}(\cat{Mod} B^{\mrm{op}})$,
restricting to a duality between 
$\msf{D}_{\mrm{f}}^{\mrm{b}}(\cat{Mod} A)$ and 
$\msf{D}_{\mrm{f}}^{\mrm{b}}(\cat{Mod} B^{\mrm{op}})$.

\begin{dfn} [\cite{Ye2}, \cite{YZ1}] \label{dfn2.2} 
Let $R$ be a dualizing complex over $(A,B)$. We 
say that $R$ has the {\em Auslander property}, or that $R$ is an 
{\em Auslander dualizing complex}, if the conditions below hold.
\begin{enumerate}
\rmitem{i} For every finite $A$-module $M$, every 
integers $p > q$, and every $B^{\mrm{op}}$-sub\-module  
$N \subset \opn{Ext}^p_A(M, R)$, one has 
$\opn{Ext}^q_{B^{\mrm{op}}}(N, R) = 0$.
\rmitem{ii} The same holds after exchanging $A$ and $B^{\mrm{op}}$.
\end{enumerate}
\end{dfn}

Rings with Auslander dualizing complexes can be viewed as a 
generalization of Auslander regular rings (cf.\ \cite{Bj} and
\cite{Le}).

\begin{exa}
If $A$ is either the $n$th Weyl algebra or the universal 
enveloping algebra of a finite dimensional Lie algebra, 
then $A$ is Auslander regular, and the bimodule $R := A$ 
is an Auslander dualizing complex. 
\end{exa}

\begin{dfn} \label{dfn2.5}
An {\em exact dimension function} on $\cat{Mod} A$ is a function
\[ \opn{dim}: \cat{Mod} A \to \{ -\infty \} \cup \mbb{R}
\cup \{ \text{infinite ordinals}\} , \]
satisfying the following axioms:
\begin{enumerate}
\rmitem{i} $\opn{dim} 0 = -\infty$.
\rmitem{ii} For every short exact sequence
$0 \to M' \to M \to M'' \to 0$ one has
$\opn{dim} M = \max\, \{ \opn{dim} M', \opn{dim} M'' \}$.
\rmitem{iii} If $M = \bigcup_{\alpha} M_{\alpha}$ then
$\opn{dim} M  = \opn{sup}\, \{ \opn{dim} M_{\alpha} \}$.
\end{enumerate}
\end{dfn} 

The basic examples of dimension functions are the Gelfand-Kirillov
dimension, denoted by $\opn{GKdim}$, and the Krull dimension, 
denoted by $\opn{Kdim}$. See \cite[Section 6.8.4]{MR}.
Here is another dimension function.

\begin{dfn}[\cite{Ye2}, \cite{YZ1}] \label{dfn2.6}
Let $R$ be an Auslander dualizing complex over $(A, B)$.
Given a finite $A$-module $M$ the {\em canonical
dimension of $M$ with respect to $R$} is
\[ \opn{Cdim}_{R; A} M := 
- \inf\, \{ q \mid \opn{Ext}^q_A(M, R) \neq 0 \}
\in {\mbb{Z}} \cup \{ -\infty \} . \]
For any $A$-module $M$ we define
\[ \opn{Cdim}_{R; A} M := \opn{sup}\, 
\{ \opn{dim} M' \mid M' \subset M \text{ is finite} \} . \]
Likewise we define $\opn{Cdim}_{R; B^{\mrm{op}}} N$ 
for a $B^{\mrm{op}}$-module $N$.
\end{dfn}

Often we shall abbreviate $\opn{Cdim}_{R; A}$ by dropping 
subscripts when no confusion can arise.
According to \cite[Theorem 2.1]{YZ1} 
$\opn{Cdim}$ is an exact dimension function on $\cat{Mod} A$
and $\cat{Mod} B^{\mrm{op}}$.

The following concept is due to Van den Bergh.
Since $R$ is a complex of $A^{\mrm{e}}$-modules, the complex
$R \otimes R$ consists of modules over
$A^{\mrm{e}} \otimes A^{\mrm{e}} \cong (A^{\mrm{e}})^{\mrm{e}}$. 
In the definition below 
$\opn{RHom}_{A^{\mrm{e}}}(A, R \otimes R)$ is computed using the 
``outside'' $A^{\mrm{e}}$-module structure of $R \otimes R$,
and the resulting complex retains the ``inside'' 
$A^{\mrm{e}}$-module structure.

\begin{dfn}[{\cite[Definition 8.1]{VdB}}] \label{dfn2.3}
Let $R$ be a dualizing complex over $A$. If there is an isomorphism
\[ \rho: R \to \opn{RHom}_{A^{\mrm{e}}}(A, R \otimes R) \]
in $\msf{D}(\cat{Mod} A^{\mrm{e}})$ then we call $(R, \rho)$, 
or just $R$, a {\em rigid} dualizing complex. 
The isomorphism $\rho$ is 
called a {\em rigidifying isomorphism}.
\end{dfn}

A rigid dualizing complex, if it exists, is unique up to 
isomorphism, by \cite[Proposition 8.2]{VdB}.

A ring homomorphism $A \to B$ is said to be {\em finite} 
if $B$ is a finite $A$-module on both sides.

\begin{dfn}\cite[Definition 3.7]{YZ1}
\label{dfn2.4}
Let $A \to B$ be a finite homomorphism of $\k$-algebras. Assume the
rigid dualizing complexes $(R_A, \rho_A)$ and $(R_B, \rho_B)$ exist. 
Let $\opn{Tr}_{B / A}: R_B \to R_A$ be a morphism in 
$\msf{D}(\cat{Mod} A^{\mrm{e}})$.
We say $\opn{Tr}_{B / A}$ is a {\em rigid trace} if it satisfies  
the following two conditions:
\begin{enumerate}
\rmitem{i} $\opn{Tr}_{B / A}$ induces isomorphisms
\[ R_B \cong \opn{RHom}_A(B, R_A) \cong 
\opn{RHom}_{A^{\mrm{op}}}(B, R_A) \]
in $\msf{D}(\cat{Mod} A^{\mrm{e}})$.
\rmitem{ii} The diagram
\[ \begin{CD}
R_B @>{\rho_B}>> \opn{RHom}_{B^{\mrm{e}}}(B, R_B \otimes R_B) \\
@V \opn{Tr} VV @VV{\opn{Tr} \otimes \opn{Tr}}V \\
R_A @>{\rho_A}>> \opn{RHom}_{A^{\mrm{e}}}(A, R_A \otimes R_A)
\end{CD} \]
in $\msf{D}(\cat{Mod} A^{\mrm{e}})$ is commutative.
\end{enumerate}
Often we shall say that
$\opn{Tr}_{B / A} : (R_B, \rho_B) \to (R_A, \rho_A)$
is a rigid trace morphism.
\end{dfn}

By \cite[Theorem 3.2]{YZ1}, a rigid trace $\opn{Tr}_{B / A}$ 
is unique (if it exists). In  
particular taking the identity map $A \to A$ and any two rigid 
dualizing complexes $(R, \rho)$ and $(R', \rho')$ 
over $A$, it follows there is a 
unique isomorphism $R \iso R'$ that is a rigid trace; see 
\cite[Corollary 3.4]{YZ1}.
Given another finite homomorphism $B \to C$ such that the rigid 
dualizing complex $(R_C, \rho_C)$ and the rigid trace
$\opn{Tr}_{C / B}$ exist, the composition
$\opn{Tr}_{C / A} := \opn{Tr}_{C / B} \circ \opn{Tr}_{B / A}$ 
is a rigid trace.

Finally we mention that by \cite[Corollary 3.6]{YZ1} 
the cohomology 
bimodules $\mrm{H}^i R_A$ of the rigid dualizing complex are 
central $\mrm{Z}(A)$-bimodules, where $\mrm{Z}(A)$ is the center 
of $A$.

\begin{exa} \label{exa2.1}
Suppose $A$ is a finitely generated commutative $\k$-algebra. 
Choose a finite homomorphism
$\k[\bsym{t}] \to A$ where 
$\k[\bsym{t}] = \k[t_1, \ldots, t_n]$ is the polynomial algebra. 
Define
$R_A := \opn{RHom}_{\k[\bsym{t}]}
(A, \Omega^n_{\k[\bsym{t}] / k}[n])$,
and consider this as an object of 
$\msf{D}^{\mrm{b}}(\cat{Mod} A^{\mrm{e}})$. 
By \cite[Proposition 5.7]{Ye4} the complex $R_A$ is an Auslander
rigid dualizing complex, and in fact it is equipped with a 
canonical rigidifying isomorphism $\rho_A$.
\end{exa}

\section{Quasi-Coherent Ringed Schemes and Localization}

In order to study geometric properties of dualizing complexes it 
is convenient to use the language of schemes and quasi coherent 
sheaves. See \cite[Theorems 0.1 and 0.2]{YZ3}.

Let $(X, \mcal{A})$ be a ringed space over $\k$. 
Thus $X$ is a topological space and $\mcal{A}$ is a sheaf of 
(possibly noncommutative) $\k$-algebras on $X$. By an
$\mcal{A}$-bimodule we mean a sheaf $\mcal{M}$ of $\k$-modules 
on $X$ together with a left $\mcal{A}$-module structure and a 
right $\mcal{A}$-module structure that commute with each other. 
In other words $\mcal{M}$ is a module over the sheaf of rings
$\mcal{A} \otimes_{\k_X} \mcal{A}^{\mrm{op}}$,
where $\k_X$ is the constant sheaf $\k$ on $X$.
An $\mcal{A}$-ring is a sheaf $\mcal{B}$ of rings on $X$
together with a ring homomorphism 
$\mcal{A} \to \mcal{B}$. 
Note that $\mcal{B}$ is an $\mcal{A}$-bimodule.

\begin{dfn} \label{dfn3.1}
Let $X$ be $\k$-scheme. An $\mcal{O}_{X}$-ring $\mcal{A}$ is 
called a {\em quasi-coherent $\mcal{O}_{X}$-ring}
if $\mcal{A}$ is a quasi-coherent 
$\mcal{O}_{X}$-module on both sides. The pair 
$(X, \mcal{A})$ is then called a {\em quasi-coherent 
ringed scheme}.
\end{dfn}

Let $(X, \mcal{A})$ be a quasi-coherent ringed scheme.
An $\mcal{A}$-module $\mcal{M}$ is called quasi-coherent if 
locally, on every sufficiently small open set $U$, it has a free 
resolution
\[ \mcal{A}|_U^{(J)} \to \mcal{A}|_U^{(I)} \to \mcal{M}|_U \to 
0 ; \]
cf.\ \cite{EGA-I}. Equivalently, $\mcal{M}$ is 
quasi-coherent as $\mcal{O}_X$-module. 
We shall denote the category of quasi-coherent 
$\mcal{A}$-modules by $\cat{QCoh} \mcal{A}$. 

\begin{prop} \label{prop3.16}
Let $(X,\mcal{A})$ be a quasi-coherent ringed scheme, let 
$U \subset X$ be an affine open set and
$A := \Gamma(U, \mcal{A})$. The functor $\Gamma(U, -)$
is an equivalence of categories
$\cat{QCoh} \mcal{A}|_U \to \cat{Mod} A$.
\end{prop} 

\begin{proof}
This is a slight generalization of \cite[Corollary 1.4.2 and
Theorem 1.5.1]{EGA-I}. See also \cite[Corollary II.5.5]{Ha}. 
\end{proof}

Given an $A$-module $M$ we shall usually denote the corresponding 
quasi-coherent $\mcal{A}|_U$-module by
$\mcal{A}|_U \otimes_{A} M$.

The following definition is due to Silver \cite[p.\ 47]{Si}.

\begin{dfn} 
\label{dfn3.3}
Let $A$ be a ring. An $A$-ring $A'$ is called a {\em 
localization} of $A$ if $A'$ is a flat $A$-module on both 
sides, and if the multiplication map 
$A' \otimes_{A} A' \to A'$ is bijective.
\end{dfn}

\begin{exa}
\label{exa3.4}
Let $A$ be a ring and $S \subset A$ a (left and right) denominator 
set. The ring of fractions $A_{S}$ of $A$ with respect to $S$ is the 
prototypical example of a localization of $A$. For reference we call
such a localization an {\em Ore localization}.
\end{exa}

We remind the reader that a denominator set $S$ 
is a multiplicatively closed subset of $A$ 
satisfying the left and right Ore conditions and the 
left and right torsion conditions (see \cite[Section 2.1]{MR}).
The left Ore condition is that for all $a \in A$ and $s \in S$ 
there exist $a' \in A$ and $s' \in S$ such that 
$a s' = a' s$. The left torsion condition is
\[ \{ a \in A \mid a s = 0 \text{ for some } s \in S \} \subset
\{ a \in A \mid s a = 0 \text{ for some } s \in S \}  . \]
The right Ore and torsion conditions for $A$ are the respective 
left conditions for $A^{\mrm{op}}$. 

Not all localizations are Ore, as we see in Example \ref{exa3.1}

Here is a list of some nice descent properties enjoyed by 
localization, that are proved  in \cite[Section 1]{Si}.

\begin{lem} \label{lem3.7}
Let $A$ be a ring and let $A'$ be a localization of $A$.
\begin{enumerate}
\item For any $A'$-module $M'$ the multiplication 
$A' \otimes_{A} M' \to M'$ is bijective.
\item Let $M'$ be an $A'$-module and $M \subset M'$ an 
$A$-submodule. Then the multiplication 
$A' \otimes_{A} M \to M'$ is injective.
\item Let $M$ be an $A$-module and 
$\phi: M \to A' \otimes_{A} M$ the homomorphism 
$\phi(m) := 1 \otimes m$. Then for any $A'$-submodule 
$N' \subset A' \otimes_{A} M$ the multiplication
$A' \otimes_{A} \phi^{-1}(N') \to N'$ is bijective.
\item In the situation of part \tup{(3)} the $A$-submodule
$\phi(M) \subset A' \otimes_{A} M$ is essential.
\item Localization of a left noetherian ring is left noetherian.
\end{enumerate}
\end{lem}

\begin{prop} \label{prop3.5}
Let $(X, \mcal{A})$ be a quasi-coherent ringed scheme, and
let $V \subset U$ be affine open sets in $X$. Define
$C := \Gamma(U, \mcal{O}_X)$, $C' := \Gamma(V, \mcal{O}_X)$,
$A := \Gamma(U, \mcal{A})$ and $A' := \Gamma(V, \mcal{A})$.
\begin{enumerate}
\item The multiplication maps 
$C' \otimes_C A \to A'$ and
$A \otimes_C C' \to A'$
are bijective. 
\item $A \to A'$ is a localization. 
\item Let $\mcal{M}$ be a quasi-coherent $\mcal{A}$-module. 
Then the multiplication map
\[ A' \otimes_{A} \Gamma(U, \mcal{M}) \to \Gamma(V, \mcal{M}) \]
is bijective.
\end{enumerate}
\end{prop}

\begin{proof} 
Define $C := \Gamma(U, \mcal{O}_{X})$ and 
$C' := \Gamma(V, \mcal{O}_{X})$. 
We first show that $C'$ is a localization of $C$, namely that
$\phi: C \to C'$ is flat and 
$\psi: C' \otimes_C C' \to C'$ is bijective. 
This can be checked locally on $V$. Choose an affine open covering 
$V = \bigcup_{i} V_i$ with
$V_i = \opn{Spec} C_{s_i}$ for suitable elements $s_i \in C$.
We note that 
$C_{s_i} \cong C'_{s_i} \cong \Gamma(V_i, \mcal{O}_{X})$ 
for all $i$. Hence restricting $\phi$ and $\psi$ to $V_i$, 
namely applying $C_{s_i} \otimes_C -$ to them, we obtain bijections.

Let $\mcal{M}$ be any quasi-coherent $\mcal{A}$-module. 
By \cite[Proposition 5.1]{Ha} multiplication
\[ C' \otimes_{C} \Gamma(U, \mcal{M}) \to \Gamma(V, \mcal{M}) \]
is bijective. 

Because $\mcal{A}$ is a quasi-coherent left and right 
$\mcal{O}_X$-module, the previous formula implies that
$C' \otimes_C A \to A'$ and $A \otimes_C C' \to A'$
are both bijective. In addition, since we now know that
\[ A' \otimes_A \Gamma(U, \mcal{M}) \cong
C' \otimes_C \Gamma(U, \mcal{M}) , \] 
we may conclude that 
\[ A' \otimes_A \Gamma(U, \mcal{M}) \to \Gamma(V, \mcal{M}) \]
is bijective.

Finally we have a 
sequence of isomorphisms, all compatible with the multiplication 
homomorphisms into $A'$:
\[ \begin{aligned}
A' \otimes_A A' & \cong A' \otimes_A (A \otimes_C C') 
\cong A' \otimes_C C' \\
& \cong (A \otimes_C C') \otimes_C C' 
\cong A \otimes_C (C' \otimes_C C') \\
& \cong A \otimes_C C' \cong A' . 
\end{aligned} \]
\end{proof}

\begin{exa} \label{exa3.1}
Let $X$ be an elliptic curve over $\mbb{C}$ and $O \in X$ the 
zero element for the group structure. Let $P \in X$ be any 
non-torsion point. Define
$U := X - \{ O \}$ and $V := X - \{ O, P \}$, which are affine open 
sets, and $C := \Gamma(U, \mcal{O}_{X})$,
$C' := \Gamma(V, \mcal{O}_{X})$. By the previous proposition 
$C \to C'$ is a localization. We claim this is not an Ore 
localization. If it were then there would be some non-invertible 
nonzero function $s \in C$ that becomes invertible in $C'$. Hence 
the divisor of $s$ on $X$ would be $(s) = n (O - P)$ for some 
positive integer $n$. In the group structure this would mean 
that $P$ is a torsion point, and this is a contradiction.
\end{exa}

\begin{dfn} \label{dfn3.8}
Let $A$ be a ring, let $M$ be an $A$-bimodule and let $A'$ be a 
localization of $A$. If the canonical homomorphisms
\[ A' \otimes_{A} M \to A' \otimes_{A} M \otimes_{A} A' \]
and
\[ M \otimes_{A} A' \to A' \otimes_{A} M \otimes_{A} A' \]
are bijective then $M$ is said to be {\em evenly localizable} to 
$A'$.
\end{dfn}

Trivially the bimodule $M := A$ is evenly localizable to $A'$. The 
next lemma is also easy and we omit its proof. 

\begin{lem}
\label{lem3.9} 
Let $M$ be an $A$-bimodule. Suppose $A'$ is a localization of $A$.
\begin{enumerate}
\item If $A$ is commutative, $M$ is a central $A$-bimodule and 
$A'$ is an $A$-algebra, then $M$ is evenly localizable to $A'$.
\item If there is a short exact sequence of $A$-bimodules 
$0 \to L \to M \to N \to 0$ with $L$ and $N$ evenly localizable 
to $A'$, then $M$ is also evenly localizable to $A'$.
\item Suppose $M \cong \lim_{i \to} M_i$ for some directed system
of $A$-bimodules $\{ M_i \}$. If each $M_i$ is evenly localizable 
to $A'$ then so is $M$. 
\end{enumerate}
\end{lem}

\begin{prop} 
\label{prop3.10}
Let $C$ be a commutative ring and let $M$ be a $C$-bimodule. 
Define $U := \opn{Spec} C$. 
The following conditions are equivalent.
\begin{enumerate}
\rmitem{i} For any multiplicatively closed subset $S \subset C$, 
with  localization $C_{S}$, $M$ is evenly localizable to $C_S$.
\rmitem{ii} $M$ is evenly localizable to 
$C' := \Gamma(V, \mcal{O}_{U})$ for every affine open set 
$V \subset U$. 
\rmitem{iii} There is a sheaf of $\mcal{O}_{U}$-bimodules 
$\mcal{M}$, quasi-coherent on both sides, with
$M \cong \Gamma(U, \mcal{M})$. Such $\mcal{M}$ is unique up to a 
unique isomorphism.
\end{enumerate}
\end{prop}

\begin{proof}
(i) $\Rightarrow$ (ii): 
Let us write
\[ \phi: C' \otimes_{C} M \to C' \otimes_{C} M \otimes_{C} C' . \]
As in the proof of Proposition \ref{prop3.5} 
we choose an affine open covering 
$V = \bigcup_{i} V_i$ with 
$V_i = \opn{Spec} C_{s_i}$ and
$C'_{s_i} \cong C_{s_i}$.
It suffices to show that the homomorphism $\phi_i$ gotten by 
applying $C'_{s_i} \otimes_{C'} -$ to $\phi$ (localizing on the 
left) is bijective for all $i$.
Using the hypothesis (i) with $S := \{ s_i^l \}_{l \in \mbb{N}}$ 
and the fact that $C_{s_i} \to C'_{s_i}$ is bijective we get
\[ C'_{s_i} \otimes_{C'} (C' \otimes_{C} M) \cong
C_{s_i} \otimes_{C}  M \otimes_{C} C_{s_i} \]
and
\[ \begin{aligned}
C'_{s_i} \otimes_{C'} (C' \otimes_{C} M \otimes_{C} C') 
& \cong
C_{s_i} \otimes_{C}  M  \otimes_{C} C_{s_i} \otimes_{C} C' \\
& \cong C_{s_i} \otimes_{C}  M \otimes_{C} C_{s_i} .
\end{aligned} \]
So $\phi_i$ is bijective.

Similarly one shows that
\[ M \otimes_{C} C' \to C' \otimes_{C} M \otimes_{C} C' \]
is bijective.

\medskip \noindent 
(ii) $\Rightarrow$ (i):
For any element $s \in S$ let $V := \opn{Spec} C_{s} \subset U$.
By assumption 
\[ C_{s} \otimes_C M \cong C_{s} \otimes_C M \otimes_C C_{s}
\cong M \otimes_C C_{s} . \]
Taking direct limit over $s \in S$ we get
\[ C_S \otimes_C M \cong C_S \otimes_C M \otimes_C C_S
\cong M \otimes_C C_S . \]

\medskip \noindent 
(ii) $\Rightarrow$ (iii): 
Let $\mcal{M} := \mcal{O}_U \otimes_C M$ be the sheafification 
of the (left) $C$-module $M$ to $U$. By definition $\mcal{M}$ is a 
quasi-coherent left $\mcal{O}_U$-module. 

Given an affine open set
$V \subset U$ write $C' := \Gamma(V, \mcal{O}_U)$. 
By Proposition \ref{prop3.5} the multiplication map
$C' \otimes_C M \to \Gamma(V, \mcal{M})$
is a bijection. Therefore 
$\Gamma(V, \mcal{M}) \cong C' \otimes_C M \otimes_C C'$.
Since $M$ is evenly localizable to $C'$
it follows that 
$M \otimes_C C' \to \Gamma(V, \mcal{M})$
is also bijective. We conclude that $\mcal{M}$ is also a 
quasi-coherent right $\mcal{O}_{U}$-module.

Regarding the uniqueness, suppose $\mcal{N}$ is another 
$\mcal{O}_{U}$-bimodule quasi-coherent on both sides such that
$\Gamma(U, \mcal{N}) \cong M$ as bimodules. For any affine open 
set $V$ as above we get an isomorphism of $C'$-bimodules
\[ \Gamma(V, \mcal{M}) \cong C' \otimes_C M \otimes_C C'
\cong \Gamma(V, \mcal{N}) \]
which is functorial in $V$. Therefore 
$\mcal{M} \cong \mcal{N}$ as $\mcal{O}_{U}$-bimodules.

\medskip \noindent 
(iii) $\Rightarrow$ (ii):
Since $\mcal{M}$ is quasi-coherent on both sides, for any 
affine open set $V = \opn{Spec} C'$ we have
\[ \Gamma(V, \mcal{M}) = C' \otimes_C M =
M \otimes_C C' , \]
so $M$ is evenly localizable to $C'$. 
\end{proof}

The relation between even localization and   
Ore localization of a ring is explained in the next theorem. 

\begin{thm} \label{thm3.11}
Let $C$ be a commutative ring, let $A$ be a $C$-ring and  
$S \subset C$ a multiplicatively closed subset. Denote by $C_{S}$ 
the ring of fractions of $C$ with respect to $S$. Then the 
following two conditions are equivalent.
\begin{enumerate}
\rmitem{i} The image $\bar{S}$ of $S$ in $A$ is a denominator set, 
with ring of fractions $A_{\bar{S}}$.
\rmitem{ii} The $C$-bimodule $A$ is evenly localizable to $C_S$.
\end{enumerate}
When these conditions hold the multiplication map
\[ C_S \otimes_C A \otimes_C C_S \to A_{\bar{S}} \]
is bijective.
\end{thm}

\begin{proof}
(i) $\Rightarrow$ (ii): Since $A_{\bar{S}}$ is
the left ring of fractions of $A$ with respect to $\bar{S}$
(see \cite[Section 2.1.3]{MR}), it follows that the homomorphism
$C_S \otimes_C A \to A_{\bar{S}}$ is bijective. 
On the other hand, since $A_{\bar{S}}$ is also the right ring of 
fractions, $A \otimes_C C_S \to A_{\bar{S}}$ is bijective. 

\medskip \noindent (ii) $\Rightarrow$ (i): 
Write 
\[ Q := C_S \otimes_C A \otimes_C C_S  \]
and
\[ \phi: A \to Q,\
\phi(a) := 1 \otimes a \otimes 1 . \]
The assumption that $A$ is evenly localizable to $C_S$ implies that
\[ \begin{aligned}
\opn{Ker}(\phi) & = \{ a \in A \mid a s = 0 \text{ for some }
s \in S \} \\
& = \{ a \in A \mid s a = 0 \text{ for some }
s \in S \} ,
\end{aligned} \]
verifying the torsion conditions. 

The even localization assumption also implies that
 given $a_1 \in A$ and $s_1 \in S$ there are $a_2 \in A$ 
 and $s_2 \in S$ such that
\[ s_1^{-1} \otimes a_1 \otimes 1 = 1 \otimes a_2 \otimes s_2^{-1}
\in Q . \]
Multiplying this equation by $s_1$ on the left and by $s_2$ on the 
right we obtain 
$1 \otimes a_1 s_2 \otimes 1 = 1 \otimes s_1 a_2 \otimes 1$. 
Therefore $\phi(s_1 a_2 - a_1 s_2) = 0$. 
Since $Q \cong A \otimes_C C_S$ there exists some $s_3 \in S$ 
such that
$(s_1 a_2 - a_1 s_2) s_3 = 0$ in $A$, i.e.\
\[ s_1 (a_2 s_3) = a_1 (s_2 s_3) \in A . \]
We have verified the right Ore condition. The left Ore condition 
is verified the same way.
\end{proof}

\begin{rem}
The theorem applies to any ring $A$ and any commutative 
multiplicatively closed subset $S \subset A$, since we can take
$C := \mbb{Z}[S] \subset A$. 
\end{rem}

We will need a geometric interpretation of Theorem 
\ref{thm3.11}. 

\begin{cor} \label{cor3.8}
Let $C$ be a commutative ring, let $U := \opn{Spec} C$ and let
$A$ be a $C$-ring. The following conditions are equivalent:
\begin{enumerate}
\rmitem{i} For every multiplicatively closed set $S \subset C$
the $C$-bimodule $A$ is evenly localizable to $C_S$.
\rmitem{ii} For every multiplicatively closed set $S \subset C$
its image $\bar{S} \subset A$ is a denominator set.
\rmitem{iii} There is a quasi-coherent $\mcal{O}_{U}$-ring 
$\mcal{A}$ such that $\Gamma(U, \mcal{A}) \cong A$ as $C$-rings.
\end{enumerate}
When these conditions hold the quasi-coherent $\mcal{O}_{U}$-ring 
$\mcal{A}$ is unique up to a unique isomorphism. 
\end{cor}

\begin{proof}
(i) and (ii) are equivalent by Theorem \ref{thm3.11}. The 
implication (iii) $\Rightarrow$ (i) is a special case of 
Proposition \ref{prop3.10}. It remains to show that 
(i) $\Rightarrow$ (iii).

By Proposition \ref{prop3.10} there is an
$\mcal{O}_{U}$-bimodule $\mcal{A}$, quasi-coherent on both sides, 
such that $A \cong \Gamma(U, \mcal{A})$ as $C$-bimodules. The 
bimodule $\mcal{A}$ is unique up to a unique isomorphism. 
Next by Theorem \ref{thm3.11}, for any $s \in C$, letting
$S := \{ s^i \}_{i \in \mbb{N}}$, the image
$\bar{S} \subset A$ is a denominator set. Therefore on
$V := \opn{Spec} C_{s}$ we have canonical isomorphisms
\[ \Gamma(V, \mcal{A}) \cong C_{s} \otimes_C A \otimes_C 
C_{s} \cong A_{\bar{S}} , \]
where $A_{\bar{S}}$ is the ring of fractions of $A$ with respect to 
$\bar{S}$. Hence $\mcal{A}$ 
has a unique structure of quasi-coherent $\mcal{O}_{X}$-ring.
\end{proof}

Here is a (somewhat artificial) example of a $C$-ring $A$ 
satisfying the conditions of Corollary \ref{cor3.8}, but the 
$C^{\mrm{e}}$-ring $A^{\mrm{e}}$ fails to satisfy them. 

\begin{exa} \label{exa4.3}
Let $C := \mbb{Q}[t]$ with $t$ a variable, and let
$U := \opn{Spec} C$. Take  
$A := \mbb{Q}(t)[a; \sigma]$, an Ore extension of the field
$\mbb{Q}(t)$, where $\sigma$ is the automorphism $\sigma(t) = -t$. 
Since every nonzero element $s \in C$ is invertible in $A$, 
the $C$-ring $A$ is evenly localizable to $C_S$ for any 
multiplicatively closed subset $S \subset C$. 
Hence there is a quasi-coherent ringed scheme $(U, \mcal{A})$
with $\Gamma(U, \mcal{A}) \cong A$ as $C$-rings. (In fact 
$\mcal{A}$ is a constant sheaf on $U$.) Likewise there's a 
quasi-coherent ringed scheme $(U, \mcal{A}^{\mrm{op}})$.

We claim that there does not exist a quasi-coherent ringed 
scheme $(U^2, \mcal{A}^{\mrm{e}})$ such that 
$\Gamma(U^2, \mcal{A}^{\mrm{e}}) \cong A^{\mrm{e}}$ as
$C^{\mrm{e}}$-rings. By Corollary \ref{cor3.8}
it suffices to exhibit a 
multiplicatively closed subset $S \subset C^{\mrm{e}}$ that is not 
a denominator set in $A^{\mrm{e}}$. Consider the element
$s := t \otimes 1 - 1 \otimes t \in C^{\mrm{e}}$ and
the set $S := \{ s^n \}_{n \in \mbb{N}}$. Let 
$\mu : A^{\mrm{e}} \to A$ be the multiplication map 
$\mu(a_1 \otimes a_2) := a_1 a_2$, which is a homomorphism of 
(left) $A^{\mrm{e}}$-modules, and denote by $I$ the left ideal
$\opn{Ker}(\mu)$. Then
$A^{\mrm{e}} \cdot s \subset I$. 
On the other hand 
$s (a \otimes 1) = t a \otimes 1 - a \otimes t$,
so
$\mu(s (a \otimes 1)) = t a - a t = 2 t a$, and by induction
$\mu(s^n (a \otimes 1)) = (2 t)^n a \neq 0$ for all $n \geq 0$. 
We conclude that 
$s^n (a \otimes 1) \notin A^{\mrm{e}} \cdot s$,
so $S$ fails to satisfy the left Ore condition.
\end{exa}

\begin{exa} \label{exa4.4}
The quasi-coherent ringed scheme $(U, \mcal{A})$ of
the previous example also has the following 
peculiarity: the $C^{\mrm{e}}$-module $A$ is not supported on the 
diagonal $\Delta(U) \subset U^2$. Indeed, for every $n \geq 0$ one 
has
$s^n a = (2 t)^n a \neq 0$. 
\end{exa}

\begin{dfn} \label{dfn4.3}
Let $C$ be a commutative $\k$-algebra and $M$ a
$C$-bimodule. A {\em differential $C$-filtration} on $M$ is an 
exhaustive, bounded below filtration $F = \{ F_{i} M \}$ where 
each $F_i M$ is a $C$-sub-bimodule, and $\opnt{gr}^{F} M$ is a 
central $C$-bimodule. If $M$ admits some differential 
$C$-filtration then we call $M$ a 
{\em differential $C$-bimodule}.
\end{dfn}

Localization of a ring was defined in Definition \ref{dfn3.3}, 
and even localization of a bimodule was introduced in 
Definition \ref{dfn3.8}.

\begin{prop} \label{prop4.2}
Let $C$ be a commutative $\k$-algebra and let $M$ be a differential 
$C$-bimodule. If $C'$ is a localization of $C$, then $M$ is evenly 
localizable to $C'$. 
\end{prop}

\begin{proof} 
If $M$ is a central $C$-bimodule then according to 
Lemma \ref{lem3.9}(1) $M$ is evenly localizable to $C'$.

Now let $M$ be a $C$-bimodule equipped with a differential 
$C$-filtration $F$. Say $F_{i_{0} - 1} M = 0$. 
We prove by induction on $i \geq i_{0}$
that $F_{i} M$ is evenly localizable to $C'$. First 
$F_{i_0} M$ is central, so the above applies to it.
For any $i$ there is an exact sequence
\[ 0 \to F_{i - 1} M \to F_{i} M \to \opnt{gr}^{F}_{i} M
\to 0 .\]
By the previous paragraph and by the induction hypothesis
$F_{i - 1} M$ and $\opnt{gr}^{F}_{i} M$ are evenly localizable 
to $C'$. The flatness of $C \to C'$ extends this to 
$F_{i} M$, see Lemma \ref{lem3.9}(2). Finally we use 
Lemma \ref{lem3.9}(3).
\end{proof}

\begin{cor} \label{cor4.2}
Let $C$ be a finitely generated commutative $\k$-algebra, 
$A$ a differential $C$-ring of finite type, $s \in C$
and $S := \{ s^i \}_{i \in \mbb{N}}$. Then:
\begin{enumerate}
\item The image $\bar{S}$ of $S$ in $A$ is a denominator set.
\item Let $C_s$ and $A_s$ be the Ore localizations w.r.t.\ $S$. 
Then $C_s$ is is a finitely generated $\k$-algebra and $A_s$ is a 
differential $C_s$-ring of finite type. 
\end{enumerate}
\end{cor}

\begin{proof}
(1) Use Proposition \ref{prop4.2} and Theorem \ref{thm3.11}.

\medskip \noindent
(2) Suppose $F = \{ F_i A \}$ is a differential $C$-filtration of 
finite type. Then setting
\[ F_i A_s := C_s \otimes_C (F_i A) \otimes_C C_s 
\subset A_s \]
we obtain a filtration $F$ of $A_s$ such that
$\opnt{gr}^F A_s \cong C_s \otimes_C \opnt{gr}^F A$
as graded $C_s$-algebras. 
\end{proof}

\begin{rem} \label{rem11.1}
The ideas in \cite[Theorem 4.9]{KL} can be used to show
the following. In the setup of the previous corollary let $M$ be 
a finite $A$-module and 
$\bar{M} := \opn{Im}(M \to A_s \otimes_A M)$. 
Then
\[ \opn{GKdim}_A \bar{M} = 
\opn{GKdim}_{A_s} (A_s \otimes_A M) . \]
\end{rem}

\begin{cor} \label{cor4.1}
Let $C$ be a commutative $\k$-algebra, let $U := \opn{Spec} C$
and let $A$ be a differential $C$-ring. 
\begin{enumerate}
\item There is a quasi-coherent $\mcal{O}_U$-ring 
$\mcal{A}$, unique up to a unique isomorphism, such that 
$\Gamma(U, \mcal{A}) \cong A$ as $C$-rings.
\item For any multiplicatively closed set $S \subset C$
its image $\bar{S} \subset A$ is a denominator set. 
\item Given an affine open set
$V \subset U$ let $C' := \Gamma(V, \mcal{O}_U)$. 
Then $A' := C' \otimes_C A \otimes_C C'$
is a $\k$-algebra and $A \to A'$ is a localization of rings. If 
$A$ is noetherian then so is $A'$.
\end{enumerate}
\end{cor}

\begin{proof}
$A$ is a differential $C$-bimodule, so by Proposition \ref{prop4.2} 
$A$ is evenly localizable to $C_S$ for any multiplicatively 
closed set $S \subset C$. By Corollary \ref{cor3.8} there is a 
quasi-coherent $\mcal{O}_U$-ring $\mcal{A}$, and by
Proposition \ref{prop3.5} we have 
$A' \cong \Gamma(V, \mcal{A})$. 
The assertion about noetherian rings follows from Lemma 
\ref{lem3.7}(5).
\end{proof}

\begin{prop} \label{prop4.1}
Let $C$ be a commutative $\k$-algebra, let 
$U := \opn{Spec} C$, let $M$ be a $C^{\mrm{e}}$-module and let
$\mcal{M} := \mcal{O}_{U^2} \otimes_{C^{\mrm{e}}} M$, the 
quasi-coherent $\mcal{O}_{U^2}$-module associated to $M$. Assume 
$C^{\mrm{e}}$ is noetherian. Then the 
following conditions are equivalent:
\begin{enumerate}
\rmitem{i} $M$ is a differential $C$-bimodule.
\rmitem{ii} $\mcal{M}$ is supported on the diagonal
$\Delta(U) \subset U^2$.
\end{enumerate}
\end{prop}

\begin{proof}
(i) $\Rightarrow$ (ii): Denote by
$I := \opn{Ker}(C^{\mrm{e}} \surj C)$ and
$\mcal{I} := \mcal{O}_{U^2} \otimes_{C^{\mrm{e}}} I$.
So $\mcal{I}$ is an ideal defining the diagonal $\Delta(U)$. 
Suppose $F = \{ F_i M \}$ is a differential $C$-filtration of $M$, 
with $F_{i_0 - 1} M = 0$. Then for all $i \geq i_0$ we have
\[ I^{i - i_0 + 1} \cdot F_i M = 0 . \]
It follows that the $\mcal{O}_{U^2}$-module
$F_i \mcal{M} := \mcal{O}_{U^2} \otimes_{C^{\mrm{e}}} F_i M$
is supported on $\Delta(U)$. But 
$\mcal{M} = \bigcup F_i \mcal{M}$. 

\medskip \noindent
(ii) $\Rightarrow$ (i): Let $\{ \mcal{M}_{\alpha} \}$ be the set 
of coherent $\mcal{O}_{U^2}$-submodules of $\mcal{M}$, 
so $\mcal{M} = \bigcup \mcal{M}_{\alpha}$. Now $\mcal{M}_{\alpha}$ 
is a coherent $\mcal{O}_{U^2}$-module supported on the diagonal 
$\Delta(U)$, so there is some integer $i_\alpha \geq 0$ such that 
$\mcal{I}^{i_\alpha + 1} \cdot \mcal{M}_{\alpha} = 0$. 
It follows that the $C^{\mrm{e}}$-module
$M_\alpha := \Gamma(U^2, \mcal{M}_{\alpha})$
satisfies $I^{i_\alpha + 1} \cdot M = 0$. And
$M = \bigcup M_\alpha$. 

Define a filtration $F$ on $M$ by
$F_i M := \opn{Hom}_{C^{\mrm{e}}}(C^{\mrm{e}} / I^{i + 1}, M)$
for $i \geq 0$, and $F_{-1} M := 0$. Then 
$M_\alpha \subset F_{i_\alpha} M$, and this implies that 
$M = \bigcup F_i M$. Finally 
$I \cdot F_i M \subset F_{i - 1} M$, and hence 
$\opnt{gr}^{F}_i M$ is a central $C$-bimodule.
\end{proof}

\section{Localization of Dualizing Complexes}

In this section we study the behavior of rigid dualizing complexes 
over rings with respect to localization (cf.\ Definition 
\ref{dfn3.8}).

\begin{dfn} \label{dfn6.1}
Let $A \to A'$ be a localization homomorphism between two noetherian
$\k$-algebras. Suppose the rigid dualizing complexes $(R, \rho)$ 
and $(R', \rho')$ of $A$ and $A'$ respectively exist. A {\em rigid 
localization morphism} is a morphism
\[ \mrm{q}_{A' / A} : R \to R' \]
in $\msf{D}(\cat{Mod} A^{\mrm{e}})$ satisfying the conditions
below. 
\begin{enumerate}
\rmitem{i} The morphisms $A' \otimes_A R \to R'$ 
and $R \otimes_A A' \to R'$ induced by $\mrm{q}_{A' / A}$
are isomorphisms.
\rmitem{ii} The diagram
\[ \begin{CD}
R @>{\rho}>> \opn{RHom}_{A^{\mrm{e}}}(A, R \otimes R) \\
@V \opn{q} VV @VV{\opn{q} \otimes \opn{q}}V \\
R' @>{\rho'}>> \opn{RHom}_{(A')^{\mrm{e}}}(A', R' \otimes R')
\end{CD} \]
in $\msf{D}(\cat{Mod} A^{\mrm{e}})$ is commutative, where
$\mrm{q} := \mrm{q}_{A' / A}$. 
\end{enumerate}
We shall sometimes express this by saying that
$\mrm{q}_{A' / A} : (R, \rho) \to (R', \rho')$
is a rigid localization morphism.
\end{dfn}

Here is a generalization of \cite[Theorem 3.8]{YZ2}.

\begin{thm} \label{thm6.2}
Let $A$ be a noetherian $\k$-algebra and let $A'$ be a 
localization of $A$. Assume $A$ has a dualizing complex $R$ 
such that the cohomology bimodules $\mrm{H}^i R$ are evenly 
localizable to $A'$. Then:
\begin{enumerate}
\item The complex
\[ R' := A' \otimes_A R \otimes_A A' \]
is a dualizing complex over $A'$.
\item If $R$ is an Auslander dualizing complex over $A$
then $R'$ is an Auslander dualizing complex over $A'$.
\item Suppose $R$ is a rigid dualizing complex over $A$ 
with rigidifying isomorphism $\rho$, and $A^{\mrm{e}}$ is 
noetherian. Then $R'$ is a rigid dualizing complex over $A'$.
Furthermore $R'$ 
has a unique rigidifying isomorphism $\rho'$ such that the 
morphism $\mrm{q}_{A' / A} : R \to R'$ defined by
$r \mapsto 1 \otimes r \otimes 1$ is a rigid localization 
morphism.
\item In the situation of part \tup{(3)} the rigid 
localization morphism
$\mrm{q}_{A' / A} : (R, \rho) \to (R', \rho')$
is unique.
\end{enumerate}
\end{thm}

\begin{proof}
(1) This follows essentially from the proof of 
\cite[Theorem 1.13]{YZ1}.
There $A$ was commutative and $A'$ the localization of $A$ at some 
prime ideal, but the arguments are valid for an arbitrary  
localization $A'$. 

\medskip \noindent (2)
To check the Auslander property for $R'$ let $M'$ be any finite 
$A'$-module. By Lemma \ref{lem3.7}(2) there is a finite 
$A$-module $M$ such that $M' \cong A' \otimes_{A} M$. 
For any $i$,  \cite[Lemma 3.7(1)]{YZ2} implies that 
\[ \opn{Ext}^{i}_{A'}(M', R') \cong
\opn{Ext}^{i}_{A}(M, R') \cong
\opn{Ext}^{i}_{A}(M, R) \otimes_{A} A' \]
as $(A')^{\mrm{op}}$-modules.
Given any $(A')^{\mrm{op}}$-submodule
$N' \subset \opn{Ext}^{i}_{A'}(M', R')$,
Lemma \ref{lem3.7}(3) tells us that there is an 
$A^{\mrm{op}}$-submodule
$N \subset \opn{Ext}^{i}_{A}(M, R)$
such that $N' \cong N \otimes_{A} A'$. For such $N$ we have
\[ \opn{Ext}^{j}_{(A')^{\mrm{op}}}(N', R') \cong
A' \otimes_{A} \opn{Ext}^{j}_{A^{\mrm{op}}}(N, R)  \]
which is $0$ for all $j < i$. By symmetry we get the other half of 
the Auslander property for $R'$.

\medskip \noindent 
(3) As in the proof of \cite[Theorem 3.8(2)]{YZ2} we have a 
canonical isomorphism
\[ A' \otimes_A \opn{RHom}_{A^{\mrm{e}}}(A, R \otimes R) 
\otimes_A A' \iso
\opn{RHom}_{(A')^{\mrm{e}}}(A', R' \otimes R') \]
in $\msf{D}(\cat{Mod}\, (A')^{\mrm{e}})$. This defines a 
rigidifying isomorphism $\rho'$ that respects $\mrm{q}_{A' / A}$ as 
depicted in the diagram in Definition \ref{dfn6.1}. 
Given any other morphism 
\[ \rho'' : R' \to
\opn{RHom}_{(A')^{\mrm{e}}}(A', R' \otimes R') \]
that renders the diagram commutative, applying the functor 
$A' \otimes_A - \otimes_A A'$ to the whole diagram we deduce that
$\rho^\flat = \rho'$.

\medskip \noindent (4)
Write $\mrm{q}_1 := \mrm{q}_{A / A'}$. 
Suppose $\mrm{q}_2 : (R, \rho) \to (R', \rho')$ is another rigid 
localization morphism. Consider the commutative diagrams
\[ \begin{CD}
R @>{\rho}>> \opn{RHom}_{A^{\mrm{e}}}(A, R \otimes R) \\
@V{\mrm{q}_i}VV @VV{\mrm{q}_i \otimes \mrm{q}_i}V \\
R' @>{\rho'}>> \opn{RHom}_{(A')^{\mrm{e}}}(A', R' \otimes R')
\end{CD} \]
in $\msf{D}(\cat{Mod} A^{\mrm{e}})$. Applying the base change 
functor 
$(A')^{\mrm{e}} \otimes_{A^{\mrm{e}}} -$ 
to these diagrams we obtain diagrams
\[ \begin{CD}
A' \otimes_A R \otimes_A A' 
@>{1 \otimes \rho \otimes 1}>> 
A' \otimes_A \opn{RHom}_{A^{\mrm{e}}}(A, R \otimes R) \otimes_A A' 
\\
@V{1 \otimes \mrm{q}_i \otimes 1}VV 
@VV{1 \otimes (\mrm{q}_i \otimes \mrm{q}_i) \otimes 1}V \\
R' @>{\rho'}>> \opn{RHom}_{(A')^{\mrm{e}}}(A', R' \otimes R')
\end{CD} \]
consisting of isomorphisms in 
$\msf{D}(\cat{Mod}\, (A')^{\mrm{e}})$;
cf.\ proof of \cite[Theorem 3.8(2)]{YZ2}.
We obtain an isomorphism $\tau : R' \to R'$ such that 
\[ 1 \otimes \mrm{q}_2 \otimes 1 = \tau \circ
(1 \otimes \mrm{q}_1 \otimes 1) : 
A' \otimes_A R \otimes_A A' \to R' . \]
But then 
$\tau : (R', \rho') \to (R', \rho')$ is a rigid trace morphism. 
By \cite[Theorem 3.2]{YZ2} $\tau$ has to be the identity. This 
implies
$1 \otimes \mrm{q}_2 \otimes 1 = 1 \otimes \mrm{q}_1 \otimes 1$
and therefore $\mrm{q}_2 = \mrm{q}_1$.
\end{proof}

The next proposition guarantees that under suitable assumptions
the rigid trace localizes.  

\begin{prop} \label{prop6.2}
Let 
\[ \begin{CD}
A @>{}>> A' \\
@V{}VV @VV{}V \\
B @>{}>> B'
\end{CD} \]
be a commutative diagram of $\k$-algebras, where the horizontal 
arrows are localizations, the vertical arrows are finite, and the 
multiplication maps
$A' \otimes_A B \to B'$ and $B \otimes_A A' \to B'$ are bijective.
Assume $A, A', A^{\mrm{e}}, B, B'$ and $B^{\mrm{e}}$ are all 
noetherian. Also assume the rigid dualizing complexes 
$(R_A, \rho_A)$ and $(R_B, \rho_B)$ exist, and so does
the rigid trace morphism
$\opn{Tr}_{B / A}: R_B \to R_A$.
By Theorem \tup{\ref{thm6.2}} the complexes 
$R_{A'} := A' \otimes_A R_A \otimes_A A'$
and 
$R_{B'} := B' \otimes_B R_B \otimes_B B'$
are rigid dualizing complexes over $A'$ and $B'$ respectively,
with induced rigidifying isomorphisms $\rho_{A'}$ and 
$\rho_{B'}$. Then the morphism
\[ \opn{Tr}_{B' / A'} := 1 \otimes \opn{Tr}_{B / A} \otimes 1 :
R_{B'} \to R_{A'} \]
is a rigid trace.
\end{prop}

\begin{proof}
We begin by showing that the morphism
$\psi': R_{B'} \to \opn{RHom}_{A'}(B', R_{A'})$
induced by $\opn{Tr}_{B' / A'}$ is an isomorphism. Let's recall 
how $\psi'$ is defined: one chooses a quasi-isomorphism 
$R_{A'} \to I'$ where $I'$ is a bounded below complex of injective 
$(A')^{\mrm{e}}$-modules. Then $\opn{Tr}_{B' / A'}$ is represented 
by an actual homomorphism of complexes 
$\tau': R_{B'} \to I'$. The formula for 
$\psi': R_{B'} \to \opn{Hom}_{A'}(B', I')$
is $\psi'(\beta')(b') = \tau'(b' \beta')$
for $\beta' \in R_{B'}$ and $b' \in B'$. 

Let $R_{A} \to I$ be a quasi-isomorphism 
where $I$ is a bounded below complex of injective 
$A^{\mrm{e}}$-modules, and let $\tau: R_{B} \to I$
be a homomorphism of complexes representing 
$\opn{Tr}_{B / A}$. We know that the homomorphism
$\psi: R_{B} \to \opn{Hom}_{A}(B, I)$ given by the formula
$\psi(\beta)(b) = \tau(b \beta)$
is a quasi-isomorphism. 

Since $R_{A'} \cong A' \otimes_A I \otimes_A A'$ 
there is a quasi-isomorphism 
$A' \otimes_A I \otimes_A A' \to I'$, 
and using it we can assume that
$\tau' = 1 \otimes \tau \otimes 1$ 
as morphisms
\[ R_{B'} = A' \otimes_A R_B \otimes_A A' \to 
A' \otimes_A I \otimes_A A' \to I' . \]
Thus we get a commutative diagram 
\[ \begin{CD}
R_B @>{\psi}>> \opn{RHom}_{A}(B, R_A) \\
@V{}VV @V{}VV \\
R_{B'} @>{\psi'}>> \opn{RHom}_{A'}(B', R_{A'})
\end{CD} \]
in $\msf{D}(\cat{Mod} A^{\mrm{e}})$.
Applying the base change $- \otimes_{A^{\mrm{e}}} (A')^{\mrm{e}}$
to the diagram we conclude that
$\psi' = 1 \otimes \psi \otimes 1$. 
So it is an isomorphism. 

By symmetry 
$R_{B'} \to \opn{RHom}_{(A')^{\mrm{op}}}(B', R_{A'})$
is also an isomorphism.

Next we have to show that the diagram
\begin{equation} \label{eqn6.1}
\begin{CD}
R_{B'} @>{\rho_{B'}}>> \opn{RHom}_{(B')^{\mrm{e}}}
(B', R_{B'} \otimes R_{B'}) \\
@V \opn{Tr} VV @VV{\opn{Tr} \otimes \opn{Tr}}V \\
R_{A'} @>{\rho_{A'}}>> \opn{RHom}_{(A')^{\mrm{e}}}
(A', R_{A'} \otimes R_{A'})
\end{CD}
\end{equation}
is commutative. This is true since (\ref{eqn6.1}) gotten by 
applying $- \otimes_{A^{\mrm{e}}} (A')^{\mrm{e}}$
to the commutative diagram 
\[ \begin{CD}
R_{B} @>{\rho_{B}}>> \opn{RHom}_{B^{\mrm{e}}}
(B, R_{B} \otimes R_{B}) \\
@V \opn{Tr} VV @VV{\opn{Tr} \otimes \opn{Tr}}V \\
R_{A} @>{\rho_{A}}>> \opn{RHom}_{A^{\mrm{e}}}
(A, R_{A} \otimes R_{A}) .
\end{CD} \]
\end{proof}

If a $\k$-algebra $A$ has an Auslander rigid dualizing complex $R$ 
then we write 
$\opn{Cdim}_A := \opn{Cdim}_{R; A}$
for this preferred dimension function. 

We finish this section with a digression from our main theme, to 
present this corollary to Theorem \ref{thm6.2}.

\begin{cor}
\label{cor6.4}
Suppose $\opn{char} \k = 0$, $A$ is the $n$th Weyl algebra 
over $\k$ and $D$ is its total ring of fractions, i.e.\ the 
$n$th Weyl division ring. Then $D[2n]$ is an Auslander rigid 
dualizing complex over $D$, and hence the canonical dimension of 
$D$ is $\opn{Cdim}_D D = 2n$. 
\end{cor}

\begin{proof}
By \cite{Ye3}, $A[2n]$ is a rigid Auslander dualizing complex 
over $A$. Now use Theorem \ref{thm6.2} with $A' := D$. 
\end{proof}

We see that unlike the Gelfand-Kirillov dimension $\opn{GKdim}$,
that cannot distinguish between the various Weyl division rings 
(since $\opn{GKdim} D = \infty$), the canonical dimension is an 
intrinsic invariant of $D$ that does recover the number $n$. 
Moreover this fact can be expressed as a ``classical'' 
formula, namely 
\[ \opn{Ext}^i_{D \otimes \mrm{D}^{\mrm{op}}}(D, D \otimes D) \cong 
\begin{cases}
D & \text{ if } i = 2 n, \\
0 & \text{ otherwise} . 
\end{cases} \]

\section{Perverse Modules and the Auslander Condition}

In this section we discuss t-structures on the derived category
$\msf{D}^{\mrm{b}}_{\mrm{f}}(\cat{Mod} A)$. 
We begin by recalling the definition of a t-structure and its 
basic properties, following \cite[Chapter X]{KS}. 

\begin{dfn} \label{dfn8.1}
Suppose $\cat{D}$ is a triangulated category and
$\cat{D}^{\leq 0}, \cat{D}^{\geq 0}$ are two full subcategories.
Let
$\cat{D}^{\leq n} := \cat{D}^{\leq 0}[-n]$
and
$\cat{D}^{\geq n} := \cat{D}^{\geq 0}[-n]$.
We say $(\cat{D}^{\leq 0}, \cat{D}^{\geq 0})$ is a
{\em t-structure on $\cat{D}$} if:
\begin{enumerate}
\rmitem{i} $\cat{D}^{\leq -1} \subset \cat{D}^{\leq 0}$
and
$\cat{D}^{\geq 1} \subset \cat{D}^{\geq 0}$.
\rmitem{ii} $\opn{Hom}_{\cat{D}}(M, N) = 0$ for
$M \in \cat{D}^{\leq 0}$ and $N \in \cat{D}^{\geq 1}$.
\rmitem{iii} For any $M \in \cat{D}$ there is a distinguished
triangle
\[ M' \to M \to M'' \to M'[1] \]
in $\cat{D}$ with $M'\in \cat{D}^{\leq 0}$ and
$M'' \in \cat{D}^{\geq 1}$.
\end{enumerate}
When these conditions are satisfied we define
the {\em heart of $\cat{D}$} to be the full subcategory
$\cat{D}^{0} := \cat{D}^{\leq 0} \cap \cat{D}^{\geq 0}$.
\end{dfn}

Given a t-structure there are truncation functors
$\tau^{\leq n} : \cat{D} \to \cat{D}^{\leq n}$
and
$\tau^{\geq n} : \cat{D} \to \cat{D}^{\geq n}$,
and functorial morphisms
$\tau^{\leq n} M \to M$,
$M \to \tau^{\geq n} M$
and
$\tau^{\geq n + 1} M \to (\tau^{\leq n} M) [1]$
such that 
\[ \tau^{\leq n} M \to M \to \tau^{\geq n + 1} M
\to (\tau^{\leq n} M) [1] \]
is a distinguished triangle in $\cat{D}$.
One shows that the heart $\cat{D}^{0}$ is an abelian category,
and the functor
\[ \mrm{H}^{0} := \tau^{\leq 0} \tau^{\geq 0} \cong
\tau^{\geq 0} \tau^{\leq 0} : \cat{D} \to \cat{D}^{0} \]
is a cohomological functor.

\begin{exa} \label{exa8.2}
Let $A$ be a left noetherian ring. The {\em standard 
t-structure} on 
$\msf{D}^{\mrm{b}}_{\mrm{f}}(\cat{Mod} A)$ is
\[ \msf{D}^{\mrm{b}}_{\mrm{f}}(\cat{Mod} A)^{\leq 0} :=
\{ M \in \msf{D}^{\mrm{b}}_{\mrm{f}}(\cat{Mod} A) \mid
\mrm{H}^{j} M = 0 \text{ for all } j > 0 \} \]
and
\[ \msf{D}^{\mrm{b}}_{\mrm{f}}(\cat{Mod} A)^{\geq 0} :=
\{ M \in \msf{D}^{\mrm{b}}_{\mrm{f}}(\cat{Mod} A) \mid
\mrm{H}^{j} M = 0 \text{ for all } j < 0 \} . \]
For a complex
\[ M = (\cdots \to M^{n} \xrightarrow{\mrm{d}^{n}} M^{n + 1} 
\to \cdots) \]
the truncations are
\[ \tau^{\leq n} M =
(\cdots \to M^{n - 2} \to M^{n - 1} \to \opn{Ker}(\mrm{d}^{n}) \to
0 \to \cdots) \]
and
\[ \tau^{\geq n} M =
(\cdots \to 0 \to \opn{Coker}(\mrm{d}^{n - 1}) \to
M^{n + 1} \to M^{n + 2} \to \cdots) . \]
The heart $\msf{D}^{\mrm{b}}_{\mrm{f}}(\cat{Mod} A)^0$
is equivalent to $\cat{Mod}_{\mrm{f}} A$.
\end{exa}

Other t-structures on 
$\msf{D}^{\mrm{b}}_{\mrm{f}}(\cat{Mod} A)$
shall be referred to as {\em perverse t-structures}, 
and the notation 
$\bigl( {}^{p}\msf{D}^{\mrm{b}}_{\mrm{f}}(\cat{Mod} A)^{\leq 0},
{}^{p}\msf{D}^{\mrm{b}}_{\mrm{f}}(\cat{Mod} A)^{\geq 0} \bigr)$
shall be used. The letter ``$p$'' stands for ``perverse'', but 
often  it will also signify a specific perversity function 
(see below). 

Now suppose for $i = 1, 2$ we are given triangulated categories
$\cat{D}_{i}$, endowed with t-structures
$(\cat{D}_{i}^{\leq 0}, \cat{D}_{i}^{\geq 0})$. An exact functor 
$F : \cat{D}_{1} \to \cat{D}_{2}$ is called {\em t-exact} if
$F(\cat{D}_{1}^{\leq 0}) \subset \cat{D}_{2}^{\leq 0}$
and
$F(\cat{D}_{1}^{\geq 0}) \subset \cat{D}_{2}^{\geq 0}$.
The functor
$F : \cat{D}_{1}^{0} \to \cat{D}_{2}^{0}$
between these abelian categories is then exact.
To apply this definition to a contravariant functor $F$ we note that
$((\cat{D}^{\geq 0})^{{\mrm{op}}}, (\cat{D}^{\leq 0})^{{\mrm{op}}})$
is a t-structure on the opposite category $\cat{D}^{{\mrm{op}}}$. 
A contravariant triangle functor 
$F : \cat{D}_{1} \to \cat{D}_{2}$ is called t-exact if
$F(\cat{D}_{1}^{\leq 0}) \subset \cat{D}_{2}^{\geq 0}$
and
$F(\cat{D}_{1}^{\geq 0}) \subset \cat{D}_{2}^{\leq 0}$.

\begin{exa} \label{exa8.3}
Let $A$ be a left noetherian $\k$-algebra and let $B$ be a right 
noetherian $\k$-algebra. Suppose we are given a dualizing complex 
$R$ over $(A, B)$, and let $\mrm{D}$ and $\mrm{D}^{{\mrm{op}}}$ 
be the duality functors $R$ induces; see Definitions \ref{dfn2.1} 
and \ref{dfn2.7}. Put on 
$\msf{D}^{\mrm{b}}_{\mrm{f}}(\cat{Mod} B^{{\mrm{op}}})$
the standard t-structure (see Example \ref{exa8.2}). Define 
subcategories
\[ {}^{p}\msf{D}^{\mrm{b}}_{\mrm{f}}(\cat{Mod} A)^{\leq 0} :=
\{ M \in \msf{D}^{\mrm{b}}_{\mrm{f}}(\cat{Mod} A) \mid
\mrm{D} M \in 
\msf{D}^{\mrm{b}}_{\mrm{f}}(\cat{Mod} B^{{\mrm{op}}})^{\geq 0}
\} \]
and
\[ {}^{p}\msf{D}^{\mrm{b}}_{\mrm{f}}(\cat{Mod} A)^{\geq 0} :=
\{ M \in \msf{D}^{\mrm{b}}_{\mrm{f}}(\cat{Mod} A) \mid
\mrm{D} M \in 
\msf{D}^{\mrm{b}}_{\mrm{f}}(\cat{Mod} B^{{\mrm{op}}})^{\leq 0}
\} . \]
Since
$\mrm{D} : \msf{D}^{\mrm{b}}_{\mrm{f}}(\cat{Mod} A) \to 
\msf{D}^{\mrm{b}}_{\mrm{f}}(\cat{Mod} B^{{\mrm{op}}})$
is a duality it follows that 
\[ \bigl (
{}^{p}\msf{D}^{\mrm{b}}_{\mrm{f}}(\cat{Mod} A)^{\leq 0},
{}^{p}\msf{D}^{\mrm{b}}_{\mrm{f}}(\cat{Mod} A)^{\geq 0}
\bigr) \]
is a t-structure on $\msf{D}^{\mrm{b}}_{\mrm{f}}(\cat{Mod} A)$,
which we call the {\em perverse t-structure induced by $R$}. 
The functors $\mrm{D}$ and $\mrm{D}^{{\mrm{op}}}$ are t-exact, and 
\[ \mrm{D} : {}^{p}\msf{D}^{\mrm{b}}_{\mrm{f}}(\cat{Mod} A)^{0} \to 
\msf{D}^{\mrm{b}}_{\mrm{f}}(\cat{Mod} B^{{\mrm{op}}})^{0} 
\approx \cat{Mod}_{\mrm{f}} B^{{\mrm{op}}} \]
is a duality of abelian categories. 
\end{exa}

\begin{dfn} \label{dfn8.3}
Suppose $A$ is a noetherian $\k$-algebra with rigid dualizing 
complex $R_A$. The perverse t-structure induced on
$\msf{D}^{\mrm{b}}_{\mrm{f}}(\cat{Mod} A)$
by $R_A$ is called the {\em rigid perverse t-structure}.
An object 
$M \in {}^{p}\msf{D}^{\mrm{b}}_{\mrm{f}}(\cat{Mod} A)^{0}$
is called a {\em perverse $A$-module}.
\end{dfn}

In the remainder of this section we concentrate on another 
method of producing t-structures on 
$\msf{D}^{\mrm{b}}_{\mrm{f}}(\cat{Mod} A)$. This method
is of a geometric nature, and closely resembles the t-structures 
that originally appeared in \cite{BBD}. 

A {\em perversity} is a function $p : \mbb{Z} \to \mbb{Z}$
satisfying
$p(i) - 1 \leq p(i + 1) \leq p(i)$.
We call the function $p(i) = 0$ the {\em trivial perversity}, 
and the function $p(i) = -i$ is called the {\em minimal perversity}. 

Let $A$ be a ring. Fix an exact dimension function 
$\opn{dim}$ on $\cat{Mod} A$ (see Definition \ref{dfn2.5}).
For an integer $i$ let $\cat{M}_{i}(\opn{dim})$ 
be the full subcategory of $\cat{Mod} A$ consisting of the modules 
$M$ with $\opn{dim} M \leq i$.
The subcategory $\cat{M}_{i}(\opn{dim})$
is localizing, and there is a functor 
\[ \Gamma_{\cat{M}_{i}(\opn{dim})} : \cat{Mod} A \to 
\cat{Mod} A \]
defined by
\[ \Gamma_{\cat{M}_{i}(\opn{dim})} M :=
\{ m \in M \mid \opn{dim} A m \leq i \} \subset M . \]
The functor $\Gamma_{\cat{M}_{i}(\opn{dim})}$ has a derived 
functor
\[ \mrm{R} \Gamma_{\cat{M}_{i}(\opn{dim})} : 
\msf{D}^+(\cat{Mod} A) \to \msf{D}^+(\cat{Mod} A) \]
calculated using injective resolutions.
For $M \in \msf{D}^+(\cat{Mod} A)$ the {\em $j$th cohomology of $M$ 
with supports in $\cat{M}_{i}(\opn{dim})$} is defined to be
\[ \mrm{H}^j_{\cat{M}_{i}(\opn{dim})} M :=
\mrm{H}^j \mrm{R} \Gamma_{\cat{M}_{i}(\opn{dim})} M . \]
The definition above was introduced in \cite{Ye2} and \cite{YZ2}. 
It is based on the following geometric paradigm. 

\begin{exa}
If $A$ is a commutative noetherian ring of finite Krull dimension
and we set $\opn{dim} M := \opn{dim} \opn{Supp} M$ 
for a finite module $M$, then
\[ \mrm{H}^j_{\cat{M}_{i}(\opn{dim})} M \cong 
\lim_{\to} \mrm{H}^j_{Z} M \]
where $Z$ runs over the closed sets in $\opn{Spec} A$ of dimension 
$\leq i$.
\end{exa}

\begin{dfn} \label{dfn8.2}
Let $A$ be a left noetherian ring. 
Given an exact dimension function $\opn{dim}$ on 
$\cat{Mod} A$ and a perversity $p$, define subcategories 
\[ {}^{p}\msf{D}^{\mrm{b}}_{\mrm{f}}(\cat{Mod} A)^{\leq 0} := 
\{ M \in \msf{D}^{\mrm{b}}_{\mrm{f}}(\cat{Mod} A) \mid
\opn{dim} \mrm{H}^{j} M < i \text{ for all } i, j \text{ with }
j > p(i) \} \]
and
\[ {}^{p}\msf{D}^{\mrm{b}}_{\mrm{f}}(\cat{Mod} A)^{\geq 0} :=
\{ M \in \msf{D}^{\mrm{b}}_{\mrm{f}}(\cat{Mod} A) \mid
\mrm{H}^{j}_{\cat{M}_{i}(\opn{dim})} M = 0
\text{ for all } i, j \text{ with } j < p(i) \} \]
of $\msf{D}^{\mrm{b}}_{\mrm{f}}(\cat{Mod} A)$. 
\end{dfn}

\begin{exa}
Suppose $\opn{dim}$ is any exact dimension function such that 
$\opn{dim} M > -\infty$ for all $M \neq 0$. Take the trivial 
perversity $p(i) = 0$. Then 
\[ {}^{p}\msf{D}^{\mrm{b}}_{\mrm{f}}(\cat{Mod} A)^{\leq 0} = 
\msf{D}^{\mrm{b}}_{\mrm{f}}(\cat{Mod} A)^{\leq 0} \]
and
\[ {}^{p}\msf{D}^{\mrm{b}}_{\mrm{f}}(\cat{Mod} A)^{\geq 0} = 
\msf{D}^{\mrm{b}}_{\mrm{f}}(\cat{Mod} A)^{\geq 0} , \]
namely the standard t-structure on 
$\msf{D}^{\mrm{b}}_{\mrm{f}}(\cat{Mod} A)$.
\end{exa}

The following lemma is straightforward.

\begin{lem} \label{lem8.3}
In the situation of Definition \tup{\ref{dfn8.2}}
let $p$ be the minimal perversity, namely $p(i) = -i$, and 
let $M \in \msf{D}^{\mrm{b}}_{\mrm{f}}(\cat{Mod} A)$.
Then:
\begin{enumerate}
\item 
$M \in {}^{p} \msf{D}^{\mrm{b}}_{\mrm{f}}(\cat{Mod} A)^{\leq 0}$
if and only if $\opn{dim} \mrm{H}^{-i} M \leq i$ for all $i$.
\item 
$M \in {}^{p} \msf{D}^{\mrm{b}}_{\mrm{f}}(\cat{Mod} A)^{\geq 0}$
if and only if
$\mrm{H}^j_{\cat{M}_{i}(\opn{dim})} M = 0$
for all $j < -i$, if and only if
$\mrm{R} \Gamma_{\cat{M}_{i}(\opn{dim})} M \in
\msf{D}^{\mrm{b}}_{\mrm{f}}(\cat{Mod} A)^{\geq -i}$
for all $i$.
\end{enumerate}
\end{lem}

Recall that if $R$ is an Auslander dualizing complex over the 
rings $(A, B)$ then the canonical dimension $\opn{Cdim}_{R}$ 
(Definition \ref{dfn2.6}) is an exact dimension function on 
$\cat{Mod} A$. 

\begin{thm} \label{thm8.6}
Let $A$ be a left noetherian $\k$-algebra and $B$ a right 
noetherian $\k$-algebra. Suppose $R$ is an Auslander dualizing 
complex over $(A,B)$. Let $\opn{dim}$ be the canonical dimension 
function $\opn{Cdim}_{R; A}$ on $\cat{Mod} A$, and let $p$ be the 
minimal perversity $p(i) = -i$. Then:
\begin{enumerate}
\item The pair
\[ \bigl( {}^p\msf{D}^{\mrm{b}}_{\mrm{f}}(\cat{Mod} A)^{\leq 0}, 
{}^p\msf{D}^{\mrm{b}}_{\mrm{f}}(\cat{Mod} A)^{\geq 0} \bigr) \]
from Definition \tup{\ref{dfn8.2}} is a t-structure on 
$\msf{D}^{\mrm{b}}_{\mrm{f}}(\cat{Mod} A)$.
\item Put on 
$\msf{D}^{\mrm{b}}_{\mrm{f}}(\cat{Mod} B^{\mrm{op}})$ 
the standard t-structure. 
Then the duality functors $\mrm{D}$ and $\mrm{D}^{{\mrm{op}}}$ 
determined by $R$ \tup{(}see Definition \tup{\ref{dfn2.7})} 
are t-exact.
\end{enumerate}
\end{thm}

\begin{proof} 
In the proof we shall use the abbreviations
$\msf{D}(A) := \msf{D}^{\mrm{b}}_{\mrm{f}}(\cat{Mod} A)$
etc.

If $M \in {}^p\msf{D}(A)^{\leq -1}$ 
then
$M[-1] \in {}^p\msf{D}(A)^{\leq 0}$. 
By Lemma \ref{lem8.3}(1) we get \linebreak
$\opn{dim} \mrm{H}^{-i} (M[-1]) \leq i$. 
Changing indices we get 
$\opn{dim} \mrm{H}^{-i} M \leq i - 1 \leq i$. 
Again using Lemma \ref{lem8.3}(1) we see that 
the first part of condition (i) of Definition 
\ref{dfn8.1} is verified. The second part of condition (i) 
is verified similarly using Lemma \ref{lem8.3}(2). 

By the Auslander condition
$\dim \opn{H}^{-j} R = \dim \opn{Ext}^{-j}_{B^{\mrm{op}}}(B, R) 
\leq j$ 
for all $j$. Therefore according to Lemma \ref{lem8.3}(1) we get
$R \in {}^p\msf{D}(A)^{\leq 0}$. 
On the other hand since
\[ \mrm{H}^{j}_{\cat{M}_{i}(\opn{dim})} R \cong
\varinjlim_{\mfrak{a} \in \mfrak{F}}
\opn{Ext}^{j}_{A}(A / \mfrak{a}, R) , \]
where $\mfrak{F}$ is the Gabriel filter of left ideals
corresponding to $\cat{M}_{i}(\opn{dim})$, the Auslander condition 
and Lemma \ref{lem8.3}(2) imply that
$R \in {}^p\msf{D}(A)^{\geq 0}$.

Let $M' \to M \to M'' \to M'[1]$ be a distinguished triangle
in $\cat{D}(A)$. Since $\opn{dim}$ is exact, and using the 
criterion in Lemma \ref{lem8.3}(1), we see that if $M'$ and $M''$
are in ${}^p\msf{D}(A)^{\leq 0}$ then so is
$M$. Likewise applying the functor 
$\mrm{H}^{j}_{\cat{M}_{i}(\opn{dim})}$
to this triangle and using Lemma \ref{lem8.3}(2) 
it follows that if $M'$ and $M''$ are in 
${}^p\msf{D}(A)^{\geq 0}$ then so is $M$.

Suppose we are given
$M \in {}^p\msf{D}(A)^{\leq 0}$,
$N \in {}^p\msf{D}(A)^{\geq 1}$
and a morphism $\phi : M \to N$. In order to prove that $\phi = 0$
we first assume $M$ is a single finite module, concentrated in some 
degree $-l$, with $l \geq 0$ and $\opn{dim} M \leq l$. 
Then $\phi$ factors through
$M \xrightarrow{\phi'} \mrm{R} \Gamma_{\cat{M}_{l}(\opn{dim})} N 
\to N$.
Now $M \in \msf{D}(A)^{\leq -l}$ and, by Lemma \ref{lem8.3}(2),
$\mrm{R} \Gamma_{\cat{M}_{l}(\opn{dim})} N \in
\msf{D}(A)^{\geq -l + 1}$;
hence $\phi' = 0$. Next let us consider the general case. Let 
$\mrm{H}^{-l} M$ be the lowest nonzero cohomology of $M$. We have 
a distinguished triangle 
\[ T := \bigl( (\mrm{H}^{-l} M)[l] \to M \to M'' \to 
(\mrm{H}^{-l} M)[l+1] \bigr) \]
where $M''$ is the standard truncation of $M$. According to 
Lemma \ref{lem8.3}(1) we have $\opn{dim} \mrm{H}^{-l} M \leq l$,
so by the previous argument the composition 
$(\mrm{H}^{-l} M)[l] \to M \xar{\phi} N$
is zero. So applying 
$\opn{Hom}_{\msf{D}(\cat{Mod} A)}(-, N)$
to the triangle $T$ we conclude that $\phi$ comes from some morphism
$\phi'': M'' \to N$. Since 
$M'' \in {}^p\msf{D}(A)^{\leq 0}$ and by induction on the 
number of nonvanishing cohomologies we have $\phi'' = 0$. 
Therefore condition (ii) is verified. 

Next suppose
$M \in \msf{D}(B^{\mrm{op}})^{\leq 0}$. 
In order to prove that
$\mrm{D}^{\mrm{op}} M \in {}^p\msf{D}(A)^{\geq 0}$
we can assume $M$ is a single finite $B^{\mrm{op}}$-module, 
concentrated in degree $-l$ for some $l \geq 0$. 
By \cite[Proposition 5.2]{YZ2} and 
its proof we deduce
\[ \begin{aligned}
\mrm{R} \Gamma_{\cat{M}_{i}(\opn{dim})} \mrm{D}^{\mrm{op}} M & =
\mrm{R} \Gamma_{\cat{M}_{i}(\opn{dim})} 
\opn{RHom}_{B^{\mrm{op}}}(M, R) \\
& \cong \opn{RHom}_{B^{\mrm{op}}}(M, 
\mrm{R} \Gamma_{\cat{M}_{i}(\opn{dim})} R) . 
\end{aligned} \]
As we saw above, the Auslander condition implies that
$\mrm{R} \Gamma_{\cat{M}_{i}(\opn{dim})} R \in
\msf{D}(B^{\mrm{op}})^{\geq -i}$. 
But 
$M \in \msf{D}(B^{\mrm{op}})^{\leq -l}$,
and hence 
\[ \opn{RHom}_{B^{\mrm{op}}}(M, 
\mrm{R} \Gamma_{\cat{M}_{i}(\opn{dim})} R) \in
\msf{D}(A)^{\geq l - i} \subset 
\msf{D}(A)^{\geq -i} . \]
Now the criterion in Lemma \ref{lem8.3}(2) tells us that 
$\mrm{D}^{\mrm{op}} M \in {}^p\msf{D}(A)^{\geq 0}$.

Let
$M \in \msf{D}(B^{\mrm{op}})^{\geq 0}$. 
We wish to prove that 
$\mrm{D}^{\mrm{op}} M \in {}^p\msf{D}(A)^{\leq 0}$.
To do so we can assume $M$ is a single finite module, 
concentrated in some degree $l \geq 0$. Then for every $i$,
$\mrm{H}^{-i} \mrm{D}^{\mrm{op}} M = 
\opn{Ext}^{-i}_{B^{\mrm{op}}}(M, R)$
has $\opn{dim} \mrm{H}^{-i} \mrm{D}^{\mrm{op}} M \leq i$.
Now apply Lemma \ref{lem8.3}(1).

At this point we know that
$\mrm{D}^{\mrm{op}}(\msf{D}(B^{\mrm{op}})^{\leq 0}) \subset
{}^p\msf{D}(A)^{\geq 0}$
and
$\mrm{D}^{\mrm{op}}(\msf{D}(B^{\mrm{op}})^{\geq 0}) \subset
{}^p\msf{D}(A)^{\leq 0}$.
Let $M \in \cat{D}(A)$
be an arbitrary complex, and consider the distinguished triangle
\[ \tau^{\leq -1} \mrm{D} M \to \mrm{D} M \to
\tau^{\geq 0} \mrm{D} M \to (\tau^{\leq -1} \mrm{D} M)[1] \]
in $\cat{D}(B^{\mrm{op}})$
gotten from the standard t-structure there.
Applying $\mrm{D}^{\mrm{op}}$ we obtain a distinguished triangle
\[ M' \to M \to M'' \to M'[1] \]
in $\cat{D}(A)$, where
$M' := \mrm{D}^{\mrm{op}} \tau^{\geq 0} \mrm{D} M$
and
$M'' := \mrm{D}^{\mrm{op}} \tau^{\leq -1} \mrm{D} M$.
This proves that condition (iii) is fulfilled, so we have a 
t-structure on $\cat{D}(A)$, and also 
the functor $\mrm{D}^{\mrm{op}}$ is t-exact. 

To finish the proof we invoke \cite[Corollary 10.1.18]{KS} 
which tells us that $\mrm{D}$ is also t-exact.
\end{proof}

\begin{que} \label{que8.8}
Let $A$ be a left noetherian ring. Find necessary and sufficient 
conditions on a dimension function $\opn{dim}$ on $\cat{Mod} A$, 
and on a perversity function $p$, such that part (1) of Theorem 
\ref{thm8.6} holds. 
\end{que}

\begin{rem} \label{rem8.2}
The idea for Definition \ref{dfn8.2} comes from 
\cite[p.\ 438 Exercise X.2]{KS}. In his recent paper \cite{Ka} 
Kashiwara considers similar t-structures. In particular his 
results imply that when $A = B$ is a commutative finitely 
generated $\k$-algebra, and $R$ is the rigid dualizing complex of 
$A$, then Theorem \ref{thm8.6} holds for any perversity 
function $p$.
(Note that here canonical dimension coincides with Krull 
dimension.) In part (2) of the theorem one has to put on
$\msf{D}^{\mrm{b}}_{\mrm{f}}(\cat{Mod} A^{\mrm{op}})$
the perverse t-structure determined by the dual perversity
$p^*(i) := -i - p(i)$. 
\end{rem}

\section{The Rigid Dualizing Complex of a Differential 
$\k$-Algebra}
\label{sec8}

We begin this section with the following consequence of previous 
work.

\begin{thm} \label{thm12.1}
Let $A$ be a differential $\k$-algebra of finite type.
Then $A$ has an Auslander rigid dualizing complex $R_{A}$. 
For a finite $A$-module $M$ the canonical dimension 
$\opn{Cdim} M$ coincides with the 
Gelfand-Kirillov dimension $\opn{GKdim} M$.
\end{thm}

\begin{proof}
According to Theorem \ref{thm11.7} $A$ has a nonnegative 
exhaustive filtration $G = \{ G_i A \}$ such that 
$\opnt{gr}^G A$ is a commutative, finitely generated, connected 
graded $\k$-algebra. Now use \cite[Corollary 6.9]{YZ1}.
\end{proof}

Recall that a ring homomorphism $f: A \to B$ is called 
{\em finite centralizing} if there exist elements 
$b_1 , \ldots, b_n \in B$ that commute with all elements of
$A$ and $B = \sum_i A b_i$. 

\begin{prop} 
\label{prop12.3}
Let $A$ be a differential $\k$-algebra of finite type and 
$f: A \to B$ a finite centralizing homomorphism. Then $B$ is 
also a differential $\k$-algebra of finite type, and the 
rigid trace $\opn{Tr}_{B / A}: R_{B} \to R_{A}$ exists.
\end{prop}

\begin{proof}
By Theorem \ref{thm11.7} we can find a differential 
$\k$-filtration of finite type $F = \{ F_{i} A \}$ of $A$
such that $\opnt{gr}^F A$ is connected. By \cite[Lemma 6.13]{YZ1} 
and its proof there is a filtration $F = \{ F_{i} B \}$ of 
$B$ such that $\opnt{gr}^F B$ is connected, 
$f(F_{i} A) \subset F_{i} B$ and
$\opnt{gr}^{F}(f): \opnt{gr}^{F} A \to \opnt{gr}^{F} B$
is a finite centralizing homomorphism. It follows that
$\opnt{gr}^{F} B$ is finite over its center, so $B$ is a  
differential $\k$-algebra of finite type. By 
\cite[Theorem 6.17]{YZ1} the 
rigid trace $\opn{Tr}_{B / A}: R_{B} \to R_{A}$ exists.
\end{proof}

Let $A$ be a differential $\k$-algebra of 
finite type with rigid dualizing complex 
$R_A$. The derived category
$\msf{D}^{\mrm{b}}_{\mrm{f}}(\cat{Mod} A)$ has on it the 
rigid perverse t-structure induced by $R_A$, whose 
heart is the category of perverse $A$-modules 
${}^{p}\msf{D}^{\mrm{b}}_{\mrm{f}}(\cat{Mod} A)^0$. 
See Definition \ref{dfn8.3}.

\begin{prop} \label{prop12.1}
Let $A \to B$ be a finite centralizing homomorphism between two 
differential $\k$-algebras of finite type. Denote by 
$\opn{rest}_{B / A} : \msf{D}(\cat{Mod} B) \to \msf{D}(\cat{Mod} A)$
the restriction of scalars functor.
\begin{enumerate}
\item Let $M \in \msf{D}^{\mrm{b}}_{\mrm{f}}(\cat{Mod} B)$. 
Then
$M \in {}^{p}\msf{D}^{\mrm{b}}_{\mrm{f}}(\cat{Mod} B)^0$
if and only if
$\opn{rest}_{B / A} M \in 
{}^{p}\msf{D}^{\mrm{b}}_{\mrm{f}}(\cat{Mod} A)^0$.
\item If $A \to B$ is surjective then the functor
\[ \opn{rest}_{B / A} : 
{}^{p}\msf{D}^{\mrm{b}}_{\mrm{f}}(\cat{Mod} B)^0
\to {}^{p}\msf{D}^{\mrm{b}}_{\mrm{f}}(\cat{Mod} A)^0 \]
is fully faithful.
\end{enumerate}
\end{prop}

\begin{proof}
(1) Define the duality functors
$\mrm{D}_A := \opn{RHom}_{A}(-, R_A)$ and 
$\mrm{D}_B :=$ \linebreak $\opn{RHom}_{B}(-, R_B)$. 
According to \cite[Proposition 3.9(1)]{YZ1} the trace 
$\opn{Tr}_{B / A} : R_B \to R_A$ gives rise to a commutative diagram
\[ \begin{CD}
\msf{D}^{\mrm{b}}_{\mrm{f}}(\cat{Mod} B)
@>{\opn{rest}_{B / A}}>> 
\msf{D}^{\mrm{b}}_{\mrm{f}}(\cat{Mod} A) \\
@VV{\mrm{D}_B}V @VV{\mrm{D}_A}V \\
\msf{D}^{\mrm{b}}_{\mrm{f}}(\cat{Mod} B^{\mrm{op}})^{\mrm{op}}
@>{\opn{rest}_{B^{\mrm{op}} / A^{\mrm{op}}}}>> 
\msf{D}^{\mrm{b}}_{\mrm{f}}(\cat{Mod} A^{\mrm{op}})^{\mrm{op}} 
\end{CD} \]
in which the vertical arrows are equivalences. 
By definition 
$M \in
{}^{p}\msf{D}^{\mrm{b}}_{\mrm{f}}(\cat{Mod} B)^0$
if and only if 
$\mrm{H}^i \mrm{D}_B M = 0$ for all $i \neq 0$.
Likewise 
$\opn{rest}_{B / A} M \in 
{}^{p}\msf{D}^{\mrm{b}}_{\mrm{f}}(\cat{Mod} A)^0$
if and only if
$\mrm{H}^i \mrm{D}_A \opn{rest}_{B / A} M = 0$ 
for all $i \neq 0$. But
\[ \mrm{H}^i \mrm{D}_A \opn{rest}_{B / A} M \cong 
\opn{rest}_{B^{\mrm{op}} / A^{\mrm{op}}} \mrm{H}^i \mrm{D}_B M . \]

\medskip \noindent
(2) In view of (1) we have a commutative diagram
\[ \begin{CD}
{}^{p}\msf{D}^{\mrm{b}}_{\mrm{f}}(\cat{Mod} B)^0
@>{\opn{rest}_{B / A}}>> 
{}^{p}\msf{D}^{\mrm{b}}_{\mrm{f}}(\cat{Mod} A)^0 \\
@VV{\mrm{D}_B}V @VV{\mrm{D}_A}V  \\
(\cat{Mod}_{\mrm{f}} B^{\mrm{op}})^{\mrm{op}}
@>{\opn{rest}_{B^{\mrm{op}} / A^{\mrm{op}}}}>> 
(\cat{Mod}_{\mrm{f}} A^{\mrm{op}})^{\mrm{op}} 
\end{CD} \]
where the vertical arrows are equivalences.
The lower horizontal arrow is a full embedding, since it 
identifies $\cat{Mod}_{\mrm{f}} B^{\mrm{op}}$ with the full 
subcategory of $\cat{Mod}_{\mrm{f}} A^{\mrm{op}}$ 
consisting of modules annihilated by 
$\opn{Ker}(A^{\mrm{op}} \to B^{\mrm{op}})$. Hence 
the top horizontal arrow is fully faithful.
\end{proof}

\begin{lem} \label{lem12.5}
Suppose $A$ and $B$ are $\k$-algebras,
$M, M' \in \msf{D}^{\mrm{b}}(\cat{Mod} A^{\mrm{e}})$
and
$N, N' \in \msf{D}^{\mrm{b}}(\cat{Mod} B^{\mrm{e}})$.
Then there is a functorial morphism
\[ \mu: \opn{RHom}_{A}(M, M') \otimes \opn{RHom}_{B}(N, N') \to
\opn{RHom}_{A \otimes B}(M \otimes N, M' \otimes N') \]
in $\msf{D}(\cat{Mod}\, (A \otimes B)^{\mrm{e}})$.
If $A$ and $B$ are left noetherian and all the modules 
$\mrm{H}^p M$ and $\mrm{H}^p N$ are finite then $\mu$ 
is an isomorphism.
\end{lem}

\begin{proof}
Choose projective resolutions
$P \to M$ and $Q \to N$ over $A^{\mrm{e}}$ and $B^{\mrm{e}}$ 
respectively. So
$P \otimes Q \to M \otimes N$ is a projective resolution over 
$(A \otimes B)^{\mrm{e}}$, and we get a map of complexes
\[ \mu: \opn{Hom}_{A}(P, M') \otimes \opn{Hom}_{B}(Q, N') \to
\opn{Hom}_{A \otimes B}(P \otimes Q, M' \otimes N') \]
over $(A \otimes B)^{\mrm{e}}$.

Now assume the finiteness of the cohomologies. To prove $\mu$ is a 
quasi-isomorphism we might as well forget the right module 
structures. Choose resolutions 
$P_{\mrm{f}} \to M$ and $Q_{\mrm{f}} \to N$
by complexes of finite projective modules over $A$ and $B$ 
respectively. We obtain a commutative diagram 
\[ \begin{CD}
\opn{Hom}_{A}(P, M') \otimes \opn{Hom}_{B}(Q, N')
@>{\mu}>> 
\opn{Hom}_{A \otimes B}(P \otimes Q, M' \otimes N') \\
@V \opn{} VV @VV{}V \\
\opn{Hom}_{A}(P_{\mrm{f}}, M') \otimes \opn{Hom}_{B}(Q_{\mrm{f}}, N')
@>{\mu}>> 
\opn{Hom}_{A \otimes B}(P_{\mrm{f}} \otimes Q_{\mrm{f}}, 
M' \otimes N') 
\end{CD} \]
in which the vertical arrows are quasi-isomorphism and the bottom 
arrow is an isomorphism of complexes. 
\end{proof}

Recall that by Proposition \ref{prop11.6} the tensor product of two 
differential $\k$-algebras of finite type is also a 
differential $\k$-algebra of finite type.  

\begin{thm} \label{thm12.3}
Suppose $A$ and $B$ are differential $\k$-algebras of finite 
type. Then the rigid dualizing complexes satisfy
\[ R_{A \otimes B} \cong R_{A} \otimes R_{B} \]
in $\msf{D}(\cat{Mod}\, (A \otimes B)^{\mrm{e}})$.
\end{thm}

\begin{proof} 
We will prove that $R_A \otimes R_B$ is a rigid dualizing complex 
over $A \otimes B$. 

Consider the Kunneth spectral sequence 
\[ (\mrm{H}^{p} R_A) \otimes (\mrm{H}^{q} R_{B}) \Rightarrow
\mrm{H}^{p + q} (R_A \otimes R_B) . \] 
Since $A \otimes B$ is noetherian it follows that 
$\mrm{H}^{p + q} (R_A \otimes R_B)$ 
is a finite $(A \otimes B)$-module on both sides.

 From Lemma \ref{lem12.5} we see that the canonical morphism
\[ A \otimes B \to
\opn{RHom}_{A \otimes B}(R_A \otimes R_B, R_A \otimes R_B) \]
in $\msf{D}(\cat{Mod}\, (A \otimes B)^{\mrm{e}})$
is an isomorphism. Likewise for 
$\opn{RHom}_{(A \otimes B)^{\mrm{e}}}$. 

Next using this lemma with $A^{\mrm{e}}$ and $B^{\mrm{e}}$, and by
the rigidity of $R_A$ and $R_B$, we get isomorphisms 
\[ \begin{aligned}
& \opn{RHom}_{(A \otimes B)^{\mrm{e}}} \bigl( A \otimes B,
(R_A \otimes R_B) \otimes (R_A \otimes R_B) \bigr)  \\
& \quad \cong 
\opn{RHom}_{A^{\mrm{e}} \otimes B^{\mrm{e}}} \bigl( A \otimes B,
(R_A \otimes R_A) \otimes (R_B \otimes R_B) \bigr)  \\ 
& \quad \cong 
\opn{RHom}_{A^{\mrm{e}}}(A, R_A \otimes R_A) \otimes 
\opn{RHom}_{B^{\mrm{e}}}(B, R_B \otimes R_B) \\
& \quad \cong R_A \otimes R_B  
\end{aligned} \]
in $\msf{D}(\cat{Mod}\, (A \otimes B)^{\mrm{e}})$.

It remains to prove that the complex $R_A \otimes R_B$ has finite 
injective dimension over $A \otimes B$ and over 
$(A \otimes B)^{\mrm{op}}$. This turns out to be quite difficult 
(cf.\ Corollary \ref{cor12.1} below).
By Theorem \ref{thm11.7} there is a filtration $F$ of $A$ 
such that $\opnt{gr}^F A$ is connected, finitely generated and 
commutative. Let $\wtil{A} := \opnt{Rees}^F A \subset A[s]$,
which is a noetherian connected graded $\k$-algebra.
By \cite[Theorem 5.13]{YZ1}, $\wtil{A}$ has a balanced
dualizing complex 
$R_{\wtil{A}} \in 
\msf{D}^{\mrm{b}}(\cat{GrMod}\, (\wtil{A})^{\mrm{e}})$.
The same holds for $B$: there's a filtration $G$, a Rees ring
$\wtil{B} := \opnt{Rees}^G B \subset B[t]$
and a balanced dualizing complex $R_{\wtil{B}}$. 
According to \cite[Theorem 7.1]{VdB} the complex
$R_{\wtil{A}} \otimes R_{\wtil{B}}$ 
is a balanced dualizing complex over $\wtil{A} \otimes \wtil{B}$. 
In particular $R_{\wtil{A}} \otimes R_{\wtil{B}}$ has finite 
graded-injective dimension over 
$\wtil{A} \otimes \wtil{B}$. 

Now $A \cong \wtil{A} / (s - 1)$, so by \cite[Lemma 6.3]{YZ1} 
the complex
\[ Q := (A \otimes \wtil{B}) \otimes_{\wtil{A} \otimes \wtil{B}} 
(R_{\wtil{A}} \otimes R_{\wtil{B}}) \]
has finite injective dimension over 
$A \otimes \wtil{B}$. But the algebra 
$A \otimes \wtil{B}$ is graded (the element $1 \otimes t$ has 
degree $1$), and $Q$ is a complex of graded 
$(A \otimes \wtil{B})$-modules. Therefore $Q$ has finite 
graded-injective dimension over this graded ring. Applying 
\cite[Lemma 6.3]{YZ1} again -- it works for any graded ring, 
connected or not -- we see that 
\[ (A \otimes B) \otimes_{A \otimes \wtil{B}} Q 
\cong (A \otimes B) \otimes_{\wtil{A} \otimes \wtil{B}} 
(R_{\wtil{A}} \otimes R_{\wtil{B}})
\]
has finite injective dimension over $A \otimes B$.

According to \cite[Theorem 6.2]{YZ1} there is an isomorphism
$R_A \cong A \otimes_{\wtil{A}} R_{\wtil{A}}[-1]$
in $\msf{D}(\cat{Mod} A)$. Likewise
$R_B \cong B \otimes_{\wtil{B}} R_{\wtil{B}}[-1]$.
Hence
\[ R_A \otimes R_B \cong 
(A \otimes B) \otimes_{\wtil{A} \otimes \wtil{B}} 
(R_{\wtil{A}} \otimes R_{\wtil{B}})[-2]  \]
has finite injective dimension over $A \otimes B$. 

By symmetry 
$R_A \otimes R_B$ has finite injective dimension also over 
$(A \otimes B)^{\mrm{op}}$.
\end{proof}

\begin{cor} \label{cor12.1}
Suppose $A$ and $B$ are differential $\k$-algebras of finite type, 
and the complexes
$M \in \msf{D}^{\mrm{b}}_{\mrm{f}}(\cat{Mod} A)$
and
$N \in \msf{D}^{\mrm{b}}_{\mrm{f}}(\cat{Mod} B)$
have finite injective dimension over $A$ and $B$ respectively. 
Then $M \otimes N$ has finite 
injective dimension over $A \otimes B$. 
\end{cor}

\begin{proof}
Let
$M^{\vee} := \opn{RHom}_A(M, R_A)$
and
$N^{\vee} := \opn{RHom}_B(N, R_B)$. 
The complexes 
$M^{\vee} \in \msf{D}^{\mrm{b}}_{\mrm{f}}(\cat{Mod} A^{\mrm{op}})$
and
$N^{\vee} \in \msf{D}^{\mrm{b}}_{\mrm{f}}(\cat{Mod} B^{\mrm{op}})$
have finite projective dimension (i.e.\ they are perfect); see
\cite[Theorem 4.5]{Ye4}. Since the tensor product of projective 
modules is projective it follows that
$M^{\vee} \otimes N^{\vee} \in
\msf{D}^{\mrm{b}}_{\mrm{f}}(\cat{Mod}\,
(A^{\mrm{op}} \otimes B^{\mrm{op}}))$
has finite projective dimension. Using Theorem \ref{thm12.3} and
Lemma \ref{lem12.5} we see that
\[ \opn{RHom}_{A \otimes B}(M \otimes N, R_{A \otimes B}) \cong
\opn{RHom}_{A \otimes B}(M \otimes N, R_{A} \otimes R_{B}) \cong
M^{\vee} \otimes N^{\vee} . \]
Applying 
$\opn{RHom}_{A \otimes B}(-, R_{A \otimes B})$
to these isomorphisms we get
\[ M \otimes N \cong 
\opn{RHom}_{A \otimes B}(M^{\vee} \otimes N^{\vee}, 
R_{A \otimes B}) , \]
so this complex has finite injective dimension.
\end{proof}

\begin{que} \label{que12.1}
Is there a direct proof of the corollary? Is it true in greater 
generality, e.g.\ for any two noetherian $\k$-algebras $A$ and 
$B$?
\end{que}

\begin{rem} \label{rem12.1}
We take this opportunity to correct a slight error in \cite{YZ1}. 
In \cite[Theorem 6.2(1)]{YZ1} the complex $R$ should be defined as
$R := (\wtil{R}_t)_0$, namely the degree $0$ component of the 
localization with respect to the element $t$. The rest of that 
theorem (including the proof!) is correct.
\end{rem}

If $A$ is a differential $\k$-algebra of finite 
type then so is the enveloping algebra $A^{\mrm{e}}$. Hence the 
rigid dualizing complex $R_{A^{\mrm{e}}}$ exists, as does
the rigid perverse t-structure on 
$\msf{D}^{\mrm{b}}_{\mrm{f}}(\cat{Mod} A^{\mrm{e}})$,
whose heart is the category 
${}^{p}\msf{D}^{\mrm{b}}_{\mrm{f}}(\cat{Mod} A^{\mrm{e}})^{0}$
of perverse $A^{\mrm{e}}$-modules.

\begin{thm} \label{thm12.2}
Let $A$ be a differential $\k$-algebra of finite type
with rigid dualizing complex $R_{A}$. Then 
$R_{A} \in 
{}^{p}\msf{D}^{\mrm{b}}_{\mrm{f}}(\cat{Mod} A^{\mrm{e}})^{0}$.
\end{thm}

\begin{proof} 
Consider the $\k$-algebra isomorphism
\[ \tau : (A^{\mrm{op}})^{\mrm{e}} = 
A^{\mrm{op}} \otimes A \iso A \otimes A^{\mrm{op}} = 
A^{\mrm{e}} \]
with formula
$\tau(a_1 \otimes a_2) := a_2 \otimes a_1$. 
Given an $A^{\mrm{e}}$-module $M$ let ${}^{\tau}M$ be the 
$(A^{\mrm{op}})^{\mrm{e}}$-module with action via $\tau$, i.e.\
\[ (a_1 \otimes a_2) \cdot_{\tau} m := 
\tau(a_1 \otimes a_2) \cdot m = a_2 m a_1 \]
for $m \in M$ and $a_1 \otimes a_2 \in (A^{\mrm{op}})^{\mrm{e}}$.
Doing this operation to the complex 
$R_A \in \msf{D}(\cat{Mod} A^{\mrm{e}})$
we obtain a complex
${}^{\tau} R_A \in \msf{D}(\cat{Mod}\, (A^{\mrm{op}})^{\mrm{e}})$.
Each of the conditions in Definitions \ref{dfn2.1} and 
\ref{dfn2.3} is automatically verified, and hence
$R_{A^{\mrm{op}}} := {}^{\tau} R_A$
is a rigid dualizing complex over $A^{\mrm{op}}$. 

According to Theorem \ref{thm12.3} we get an isomorphism
\[ R_{A^{\mrm{e}}} \cong R_A \otimes R_{A^{\mrm{op}}} =
R_A \otimes ({}^{\tau} R_A) \]
in $\msf{D}(\cat{Mod}\, (A^{\mrm{e}})^{\mrm{e}})$.
But the left (resp.\ right) $A^{\mrm{e}}$ action on 
$R_A \otimes ({}^{\tau} R_A)$ is exactly the outside 
(resp.\ inside) action on $R_A \otimes R_A$. By rigidity
(cf.\ Definition \ref{dfn2.3}) we have
\[ R_A \cong 
\opn{RHom}_{A^{\mrm{e}}}(A, R_A \otimes R_A) \cong
\opn{RHom}_{A^{\mrm{e}}}(A, R_{A^{\mrm{e}}}) \]
in $\msf{D}(\cat{Mod}\, A^{\mrm{e}})$.

Finally since $A^{\mrm{e}} \cong (A^{\mrm{e}})^{\mrm{op}}$, via 
the involution $\tau$, we may view 
$\opn{RHom}_{A^{\mrm{e}}}(-, R_{A^{\mrm{e}}})$
as an auto-duality of
$\msf{D}^{\mrm{b}}_{\mrm{f}}(\cat{Mod}\, A^{\mrm{e}})$.
By definition of the rigid t-structure this duality exchanges
$\cat{Mod}_{\mrm{f}} A^{\mrm{e}}$ and
${}^p \msf{D}^{\mrm{b}}_{\mrm{f}}(\cat{Mod}\, A^{\mrm{e}})^0$.
Since 
$A \in \cat{Mod}_{\mrm{f}} A^{\mrm{e}}$
it follows that 
$R_A \in$ \linebreak
${}^p \msf{D}^{\mrm{b}}_{\mrm{f}} (\cat{Mod}\, A^{\mrm{e}})^0$.
\end{proof}

We know that the cohomology bimodules $\mrm{H}^i R_A$ are 
central $\mrm{Z}(A)$-modules. The next lemma shall be used a few 
times. 

\begin{lem} 
\label{lem12.7}
Let $A$ be a differential $\k$-algebra of finite type and
$a \in A$ a non-invertible central regular element.
Define $B := A / (a)$. Let $R_{A}$ and $R_{B}$ denote the 
rigid dualizing complexes of $A$ and $B$ respectively. 
Then there is a long exact sequence
\[ \cdots \to \mrm{H}^{i} R_{B} \to
\mrm{H}^{i} R_{A} \xrightarrow{a}
\mrm{H}^{i} R_{A} \to \mrm{H}^{i + 1} R_{B} \to \cdots  \]
of $A$-bimodules.
\end{lem}

\begin{proof} 
Trivially $A \to B$ is a finite centralizing homomorphism.
By Proposition \ref{prop12.3} the trace morphism 
$\opn{Tr}_{B/A}: R_B \to R_A$ exists. 
In particular, $R_B \cong \opn{RHom}_A(B, R_A)$. 
There is an exact sequence of bimodules
\[ 0 \to A \xrightarrow{a} A \to B \to 0 . \]
Applying the functor $\opn{RHom}_A(-, R_A)$ to it and taking 
cohomologies we obtain the long exact sequence we want.
\end{proof}

Here are a couple of examples of differential $\k$-algebras of 
finite type and their rigid dualizing complexes.

\begin{exa}
Let $C$ be a smooth $n$-dimensional $\k$-algebra  
in characteristic $0$ and $A := \mcal{D}(C)$ the 
ring of differential operators. Then the rigid dualizing complex is
$R_{A} = A[2 n]$; see \cite{Ye3}.
\end{exa}

\begin{exa} \label{exa12.1}
Let $\mfrak{g}$ be an $n$-dimensional Lie algebra over
$\k$ and $A := \mrm{U}(\mfrak{g})$ its universal enveloping algebra.
By \cite{Ye3} the rigid dualizing complex is
$R_{A} = A \otimes (\bwedge^n \mfrak{g})[n]$,
where $\bwedge^n \mfrak{g}$ has the adjoint $A$ action on the left 
and the trivial action on the right. 
\end{exa}

Suppose $A$ is a ring with nonnegative exhaustive filtration $F$
such that the Rees ring $\wtil{A} := \opnt{Rees}^F A$ is left
noetherian. We remind that a filtered $(A, F)$-module 
$(M, F)$ is called {\em good} if it is bounded below, exhaustive, 
and $\opnt{Rees}^F M$ is a finite $\wtil{A}$-module. 

In the two previous examples the cohomology bimodules 
$\mrm{H}^i R_A$ all came equipped with filtrations that were both 
differential and good on both sides. These properties turn out to 
hold in general, as Theorems \ref{thm12.8} and \ref{thm12.9}
show.

\begin{thm} \label{thm12.8}
Let $A$ be a differential $\k$-algebra of finite type, and let 
$R_{A}$ be the rigid dualizing complex of $A$. Let $F$ be some
differential $\k$-filtration of finite type of $A$. Then 
for every $i$ there is an induced filtration $F$ of 
$\mrm{H}^{i} R_{A}$, such that $(\mrm{H}^{i} R_{A}, F)$ is a good 
filtered $(A, F)$-module on both sides.
\end{thm}

\begin{proof}
Define $\wtil{A} := \opnt{Rees}^{F} A \subset A[t]$.  
Let $\wtil{F} = \{ \wtil{F}_i \wtil{A} \}$ be the filtration from 
Lemma \ref{lem10.7}. Then
\[ \opnt{gr}^{\wtil{F}} \wtil{A} \cong (\opnt{gr}^{F} A) 
\otimes \k[t] \]
as $\k$-algebras. The center is
\[ \mrm{Z}(\opnt{gr}^{\wtil{F}} \wtil{A}) \cong 
\mrm{Z}(\opnt{gr}^{F} A) \otimes \k[t] , \]
which is a finitely generated commutative $\k$-algebra. 
Also $\opnt{gr}^{\wtil{F}} \wtil{A}$ is a finite 
$\mrm{Z}(\opnt{gr}^{\wtil{F}} \wtil{A})$-module. 
We conclude that $\wtil{F}$ is a differential $\k$-filtration
of finite type on $\wtil{A}$. Moreover each $\k$-submodule 
$\wtil{F}_i \wtil{A}$ is graded, where $\wtil{A} \subset A[t]$ 
has the grading $F$ in which $\opn{deg}^F(t) = 1$. 

Applying Theorem \ref{thm11.7} and Proposition \ref{prop11.2} to 
the filtered $\k$-algebra $(\wtil{A}, \wtil{F})$
we obtain another filtration $\wtil{G}$ on $\wtil{A}$. This new 
filtration is also differential $\k$-filtration
of finite type, and each $\k$-submodule 
$\wtil{G}_i \wtil{A} \subset \wtil{A}$ is graded (for the grading 
$F$). Furthermore $\opnt{gr}^{\wtil{G}} \wtil{A}$
is a connected graded $\k$-algebra (when considered as 
$\mbb{Z}$-graded ring with the grading $\wtil{G}$).

Define
\[ B := \opnt{Rees}^{\wtil{G}} \wtil{A} \subset \wtil{A}[s]
\subset A[s, t] . \]
This is a $\mbb{Z}^2$-graded ring with grading 
$(\wtil{G}, F)$, in which
$\opn{deg}^{(\wtil{G}, F)}(s) = (1, 0)$ and
$\opn{deg}^{(\wtil{G}, F)}(t) = (0, 1)$.

Consider the $\k$-algebra $B$ with its grading $\wtil{G}$. 
This is a connected graded $\k$-algebra. 
The quotient 
$B / (s) \cong \opnt{gr}^{\wtil{G}} \wtil{A}$
is a finitely generated commutative $\k$-algebra.
Therefore by \cite[Theorem 5.13]{YZ1},
$B$ has a balanced dualizing complex 
$R_{B} \in \msf{D}^{\mrm{b}}(\cat{GrMod} B^{\mrm{e}})$.
By \cite[Theorem 6.3]{VdB} we get
\[ R_{B} \cong 
(\opn{R} \Gamma_{\mfrak{m}} B)^{*}  \]
in $\msf{D}(\cat{GrMod} B^{\mrm{e}})$. 
Here $\Gamma_{\mfrak{m}}$ is the torsion functor 
with respect to the augmentation ideal $\mfrak{m}$ of $B$,
and 
\[ (M)^* := \opn{Hom}^{\mrm{gr}}_{\k}(M, \k) =
\boplus_i \opn{Hom}_{\k}(M_{-i}, \k) , \]
the graded dual of the graded $\k$-module $M$. 
In particular, for every $p$ there is an isomorphism 
of $B$-bimodules
$\mrm{H}^{-p} R_{B} \cong
\bigl( \mrm{H}^p_{\mfrak{m}} B \bigr)^{*}$
where
\[ \mrm{H}^p_{\mfrak{m}} B := \mrm{H}^p \mrm{R} \Gamma_{\mfrak{m}} B
\cong {\lim_{k \to}} \opn{Ext}^{p}_{B}(B / \mfrak{m}^{k}, B) . \]
Now for each $k$ 
\[ \opn{Ext}^{p}_{B}(B / \mfrak{m}^{k}, B) = 
\boplus_{(i, j) \in \mbb{Z}^2} 
\opn{Ext}^{p}_{B}(B / \mfrak{m}^{k}, B)_{(i, j)} \]
where $(i, j)$ is the $(\wtil{G}, F)$ degree. Therefore in the 
direct limit we get a double grading
\[ \mrm{H}^p_{\mfrak{m}} B = \boplus_{(i, j) \in \mbb{Z}^2} 
(\mrm{H}^p_{\mfrak{m}} B)_{(i, j)} . \]
Since for every $i$ the $\k$-module
$(\mrm{H}^p_{\mfrak{m}} B)_{i} = \boplus_{j \in \mbb{Z}} 
(\mrm{H}^p_{\mfrak{m}} B)_{(i, j)}$
is finite it follows that the graded dual 
$(\mrm{H}^p_{\mfrak{m}} B)^*$, 
which is computed with respect to the $\wtil{G}$ grading, is also 
$\mbb{Z}^2$-graded. We see that $\mrm{H}^{-p} R_{B}$ is in fact a 
$\mbb{Z}^2$-graded $B$-bimodule.

By \cite[Theorem 6.2]{YZ1} the complex
$R_{\wtil{A}} := ((R_B[-1])_s)_0$ is a rigid dualizing complex 
over the ring $\wtil{A} \cong B / (s - 1)$; cf.\ Remark 
\ref{rem12.1}.
Hence each cohomology 
\[ \mrm{H}^p R_{\wtil{A}} \cong 
\frac{\mrm{H}^{p - 1} R_{B}}{(s - 1) \cdot \mrm{H}^{p - 1} R_{B}}
= \opn{sp}^{\wtil{G}}_1 \mrm{H}^{p - 1} R_{B} \]
is a graded $\wtil{A}$-bimodule (with the $\mbb{Z}$-grading 
$F$ in which $\opn{deg}^F(t) = 1$), finite on both sides. 

Next we have $A \cong \wtil{A} / (t - 1)$. Because $t - 1$ is a 
central regular non-invertible element, Lemma \ref{lem12.7} says 
there is an exact sequence
\[ \cdots \to \mrm{H}^{p} R_{A} \to
\mrm{H}^{p} R_{\wtil{A}} \xrightarrow{t - 1}
\mrm{H}^{p} R_{\wtil{A}} \to \mrm{H}^{p + 1} R_{A} \to 
\mrm{H}^{p + 1} R_{\wtil{A}} \to \cdots  \]
of $A$-bimodules. Since $\mrm{H}^{p} R_{\wtil{A}}$ is graded the 
element $t - 1$ is a non-zero-divisor on it, and therefore 
we get an exact sequence
\[ 0 \to \mrm{H}^{p - 1} R_{\wtil{A}} \xrightarrow{t - 1}
\mrm{H}^{p - 1} R_{\wtil{A}} \to \mrm{H}^{p} R_{A} \to 0 . \]
Thus the bimodule
\[ \mrm{H}^p R_A \cong 
\frac{\mrm{H}^{p - 1} R_{\wtil{A}}}{(t - 1) \cdot 
\mrm{H}^{p - 1} R_{\wtil{A}}} =
\opn{sp}^F_1 \mrm{H}^{p - 1} R_{\wtil{A}} \]
inherits a bounded below exhaustive filtration $F$, and 
$\opn{Rees}\, (\mrm{H}^p R_A, F)$, being a quotient of 
$\mrm{H}^{p - 1} R_{\wtil{A}}$, is a finite $\wtil{A}$-module on 
both sides. By definition $(\mrm{H}^p R_A, F)$ is then
is a good filtered $(A, F)$-module on both sides.  
\end{proof}

\begin{thm} \label{thm12.9}
Let $C$ be a finitely generated commutative $\k$-algebra, 
let $A$ be differential $C$-ring of finite type, and let 
$R_{A}$ be the rigid dualizing complex of $A$. Then for every
$i$ the $C$-bimodule $\mrm{H}^{i} R_{A}$ is differential.
\end{thm}

\begin{proof} 
Using the same setup as in the proof of Theorem \ref{thm12.8}, 
define
\[ \bar{A} := \opnt{gr}^{F} A = \opn{sp}^F_{0} \wtil{A} = 
\wtil{A} / (t) . \]
So $\bar{A}$ is a $C$-algebra. Let 
$R_{\bar{A}}$ be the rigid dualizing complex of $\bar{A}$. 
By \cite[Corollary 3.6]{YZ1} 
the $\bar{A}$-bimodule $\mrm{H}^{i} R_{\bar{A}}$ 
is $Z(\bar{A})$-central, and hence it is a central $C$-bimodule.
According to Lemma \ref{lem12.7} there is an exact sequence of 
$\wtil{A}$-bimodules
\[ \mrm{H}^{i - 1} R_{\wtil{A}} \xrightarrow{t} 
\mrm{H}^{i - 1} R_{\wtil{A}} \to \mrm{H}^{i} R_{\bar{A}} . \]
Therefore 
\[ \opn{sp}_{0}\, (\mrm{H}^{i - 1} R_{\wtil{A}}) =
(\mrm{H}^{i - 1} R_{\wtil{A}}) / 
t \cdot (\mrm{H}^{i - 1} R_{\wtil{A}}) 
\inj \mrm{H}^{i} R_{\bar{A}} \]
is a central $C$-bimodule. 

To conclude the proof, consider the filtration $F$ of 
$\mrm{H}^{i} R_{A}$ from Theorem \ref{thm12.8}. Because 
$(\mrm{H}^{i} R_{A}, F) \cong 
\opn{sp}_{1}\, (\mrm{H}^{i - 1} R_{\wtil{A}})$ 
is a good filtered $(A, F)$-module, say on the left, we see that 
$(\mrm{H}^{i} R_{A}, F)$ 
is exhaustive and bounded below. Now
\[ \opnt{Rees}\, (\mrm{H}^{i} R_{A}, F) \cong 
(\mrm{H}^{i - 1} R_{\wtil{A}}) / \{ t \text{-torsion} \} , \]
so $\opnt{gr}\, (\mrm{H}^{i} R_{A}, F)$
is a quotient of 
$\opn{sp}_{0}\, (\mrm{H}^{i - 1} R_{\wtil{A}})$.
It follows that $\opnt{gr}\, (\mrm{H}^{i} R_{A}, F)$ is a 
central $C$-bimodule. Thus $F$ is a differential $C$-filtration of 
$\mrm{H}^{i} R_{A}$.
\end{proof}

\begin{cor} \label{cor8.10}
In the situation of Theorem \tup{\ref{thm12.9}} let
$U := \opn{Spec} C$. Given an affine 
open set $V \subset U$ let $C' := \Gamma(V, \mcal{O}_U)$. Then
$A' := C' \otimes_C A \otimes_C C'$ is a noetherian $\k$-algebra, 
$A \to A'$ is flat and 
$R' := A' \otimes_A R_A \otimes_A A'$ is an Auslander rigid 
dualizing complex over $A'$.
\end{cor}

\begin{proof}
According to Corollary \ref{cor4.1} $A'$ is a noetherian 
$\k$-algebra and $A \to A'$ is a localization. We know that 
$A^{\mrm{e}}$ is noetherian. By Theorem \ref{thm12.9} each of the 
cohomology bimodules $\mrm{H}^{i} R_{A}$ are differential as 
$C$-bimodules, hence by Proposition \ref{prop4.2} they are evenly 
localizable to $C'$. 
Thus all the hypothesis of Theorem \ref{thm6.2}(1,2,3) 
are satisfied.
\end{proof}

\end{document}